\documentclass[a4paper, 11pt]{article}
\usepackage{amsmath,url}
\usepackage{graphicx}
\usepackage{latexsym}
\usepackage{amssymb} % or amsfonts
\usepackage{scalerel,stackengine}
\usepackage{xcolor}
\stackMath
\newcommand\reallywidehat[1]{%
\savestack{\tmpbox}{\stretchto{%
  \scaleto{%
    \scalerel*[\widthof{\ensuremath{#1}}]{\kern-.6pt\bigwedge\kern-.6pt}%
    {\rule[-\textheight/2]{1ex}{\textheight}}%WIDTH-LIMITED BIG WEDGE
  }{\textheight}%
}{0.5ex}}%
\stackon[1pt]{#1}{\tmpbox}%
}
\renewcommand{\Bbb}{\mathbb}

\def\b{\begin{eqnarray}}
\def\e{\end{eqnarray}}

\newtheorem{theorem}{Theorem}

\begin{document}
\title{Nonlinear two-dimensional water waves with arbitrary vorticity}
\author{Delia Ionescu-Kruse\footnote{Simion Stoilow Institute of Mathematics of the Romanian Academy,
Research Unit No. 6,
 P.O. Box 1-764, RO-014700, Bucharest, Romania  and Erwin Schr\"odinger International Institute for Mathematics and Physics (ESI),
University of Vienna, Boltzmanngasse 9, 1090 Vienna, Austria;
e-mail: Delia.Ionescu@imar.ro},
        Rossen Ivanov\footnote{School of Mathematics and Statistics, Technological University Dublin, Grangegorman Lower,
Dublin D07 ADY7, Ireland and Erwin Schr\"odinger International Institute for Mathematics and Physics (ESI),
University of Vienna, Boltzmanngasse 9, 1090 Vienna, Austria;
e-mail: rossen.ivanov@tudublin.ie}
}
\maketitle

\begin{abstract}

We consider the two-dimensional water-wave  problem with a general non-zero vorticity field in a fluid volume with a flat bed and a free surface.
The nonlinear equations of motion for the chosen surface and volume variables are expressed with the aid of the Dirichlet-Neumann operator and the Green function of the Laplace operator in the fluid domain. Moreover, we provide new explicit expressions for both objects.
The field of a point vortex and its interaction with the free surface is studied as an example. In the small-amplitude long-wave Boussinesq and KdV regimes, we obtain appropriate  systems of coupled equations for the dynamics of the point vortex and the time evolution of the free surface variables.

\end{abstract}

\section{Introduction}

The study of rotational  flows of an inviscid incompressible fluid with free boundary is of high theoretical and practical interest (e. g., the evolution of surface waves on the ocean, their approximation by model equations and by computer simulations, the dynamics of wave interactions).
The  vorticity is the key quantity in the analysis of the fluid motion.
The description of the dynamics of rotational ideal fluid flows with a free boundary in two dimensions is of fundamental significance. In this paper we concentrate on the governing equations (expressed in terms of the appropriate surface and volume variables).
We provide a general description of the intricate interaction between the rotational motion in the bulk of the fluid and the surface motion. The Dirichlet-Neumann operator and the Green function of the Laplace operator in the bulk of the fluid play a very important role in the mathematical formulation of the nonlinear system of equations. We note that the Green function in this case depends on the moving boundaries of the domain, that is, the free surface of the fluid.

The problem of the Hamiltonian formulation of the hydrodynamic equations plays an important role in the water-wave theory and has a long history. A significant theoretical development has been made by Zakharov in his 1968 paper \cite{Zakharov68}. He reformulated the problem of irrotational  waves in deep water as a Hamiltonian system where the velocity potential for the irrotational flow and the surface elevation are the canonically conjugate variables that constitute the dynamics: in fact, it suffices to know the boundary values of the velocity potential at any particular time and the shape of the free surface in order to recover the flow throughout the fluid region at that time. The Hamiltonian functional can be conveniently expressed through the surface variables with the help of the non-local Dirichlet-Neumann operator, which takes boundary values of a harmonic function and returns its normal derivative.  The Dirichlet-Neumann operator also allows for a straightforward derivation of various long-wave approximations of the water-wave problem, see Craig and Groves  \cite{Craig&Groves}. Irrotational flows with free surfaces have been studied from a Hamiltonian point of view also by Miles \cite{miles} (see also Milder \cite{milder}), Benjamin and Olver  \cite{benjamin&olver} and many others.
Zakharov's  ideas  for irrotational waves in deep water have been generalized in various directions by several authors. Wahl{\'{e}}n \cite{Wahlen} provided a Hamiltonian formulation for surface waves of finite-depth water with constant vorticity using the {\it nearly}-Hamiltonian formulation of Constantin, Ivanov and Prodanov \cite{CIP}.  The formulation also involves the Dirichlet–Neumann operator.

\noindent The two-dimensional two-layer irrotational gravity water flows with a free surface and interface possesses a Hamiltonian formulation too, this has been shown by Craig, Guyenne and Kalisch \cite{CGK} (using the results from \cite{BB} and  \cite{Craig&Groves}).  For two layers with constant vorticity in each layer, the Hamiltonian approach has been extended by Constantin, Ivanov and Martin \cite{CIM}.
 For stratified  two dimensional periodic water flows with piecewise constant vorticity, accounting for Coriolis effects in the equatorial $f$-plane approximation does not hinder the Hamiltonian description of the governing equations, see Ionescu-Kruse and Martin \cite{IK-M}, Constantin and Ivanov \cite{CI}.
We  mention that  additional complexity is introduced when water-waves are propagating over a non-flat bottom (see, for example, Craig et al. \cite{CGNS},
Compelli, Ivanov and Todorov \cite{CIT},  Compelli et al. \cite{CIMT}, and the references therein).

\noindent  In  the case of an arbitrary, rotational flow, Lewis et al. have introduced in  \cite{LMMR}  non-canonical Poisson brackets and they have shown that the equations of motion for incompressible flows with a free boundary having surface tension are Hamiltonian relative to this structure.
The method used in \cite{LMMR} for obtaining the Poisson bracket is to pass from canonical brackets in the Lagrangian (material) representation to the
non-canonical brackets in the Eulerian (spatial) representation.
The non-canonical variables in \cite{LMMR} are the velocity vector field and the free surface, the Hamiltonian functional is given by the total energy.
However, this formulation is in terms of a degenerate Poisson bracket which greatly complicates the diagonalization of the bracket and the introduction of canonical variables.
In addition, in the infinite-dimensional case the Poisson bracket is defined usually only for a specific class of functionals, often called admissible
functionals. The fact that one picks out a special class of functionals is related to the fact that usually one does not provide the infinite dimensional phase space with the structure of a smooth manifold and therefore one cannot just speak of $C^\infty$-functions on this phase space. The choice of this class of functions and the definition of their functional derivatives is very subtle matter and there is no simple recipe describing how to choose these objects in our infinite dimensional system with its free boundaries. If the infinite dimensional configuration space is the diffeomorphism group of a compact oriented smooth $n$-dimensional manifold with smooth boundary, a situation which appears in the mathematical modelling of  an ideal fluid which completely fills a vessel with smooth boundary, Ebin and Marsden constructed rigorously a structure of
manifold for this configuration space in \cite{Ebin&Marsden}.  Configuration manifolds for ideal fluids with free boundaries have been introduced at a formal level in \cite{LMMR}. The Poisson bracket is defined only for the admissible functionals. In order to verify the Jacobi identity on all triples of admissible functionals, one has first to check that the bracket of two admissible functionals is again an admissible functional.
This difficulty of non-closure of the class of admissible functionals under the operation of taking Poisson brackets, was first considered by Soloviev in a series of papers (see, for example, \cite{soloviev 1993}, \cite{soloviev 2002}).

\noindent In 1859 Clebsch derived the equations of motion in an inviscid, incompressible fluid by  a variational principle. The variables introduced by Clebsch,  which appear from his original representation of the fluid velocity field and  called by some authors  Clebsch potentials (for a summary of the basic facts from the two important Clebsch papers \cite{clebsch}, see for example,  Lamb \cite{lamb}, \S 167), have the remarkable property of being canonical variables, they enable one to represent the equations of motion for ideal hydrodynamics in a Hamiltonian form. Having canonical variables one can easily calculate the Poisson brackets between any physical quantities and one also succeeds in writing down a variational principle. The choice of Clebsch variables is not unique and it is also  difficult to assign a physical meaning to them (see, for example, the discussion in Benjamin  \cite{benjamin}, \S 7 and the references therein). A variational formulation of this type has been put forward by Cotter and Bokhove in \cite{CoBo}.

Due to the difficulties outlined above, in this paper we concentrate on the equations for the rotational fluid with a free surface, rather than on the Hamiltonian description of the problem.
After the preliminaries and notations in Section 2, we reformulate the nonlinear equations of motion with the aid of the Dirichlet-Neumann operator and the Green function of the Laplace operator in the bulk of the fluid in Section 3.
If the flow is irrotational, then the velocity can be expressed as the gradient of a scalar field, namely the velocity potential. According to Weyl-Hodge theory (see \cite{Ebin&Marsden} for a summary and references), the velocity field decomposes uniquely as the sum of the gradient of the velocity potential $\varphi$ and a divergence free vector field which is tangent to the free surface. We denote by $r$ the stream function of this divergence free vector from the unique Weyl-Hodge decomposition and by $\xi$ the restriction of $\varphi$ to the free surface $\eta$.
For a nonconstant vorticity $\omega(x,z,t)$, evolving according to the vorticity equation, we recast the free boundary water-wave problem in terms of the variables $\xi$, $r$, $\eta$ and $\omega$.
The evolution in time of these variables is described by equations, involving the Dirichlet-Neumann operator and the Green function  of the Laplace operator in the bulk of the fluid. For both structures we provide explicit expressions, \eqref{G1}, \eqref{61}, respectively. As an illustration of the general setting, the small-amplitude  long-wave scaling regime is studied in Section 5, and the corresponding expression for the Dirichlet-Neumann operator is obtained \eqref{G3} and compared with the Craig-Sulem \cite{CS} expression as a power series expansion in terms of the free surface variables.
In the next Sections, the field of a point vortex and its interaction with the surface is studied as an example. In small-amplitude, long-wave regimes, the appropriate Green function for the point-vortex problem is the Green function on an infinite strip.  Alternative expressions of the Green function on an infinite strip, constructed by different methods - the method of images, the method of conformal mapping or the method of eigenfunction expression - are available in standard books on partial differential equations, see, for example,  Marchioro and Pulvirenti \cite{M&Pulvirenti}, Yu.A. Melnikov and M.Y.Melnikov \cite{Melnikov} and the references therein. We obtain an alternative form of the Green function on the strip, that is, \eqref{102}.
Finally, in the Boussinesq and KdV propagation regimes, by assuming that the point vortex  does not move under the action of its own field, we arrive at  the systems \eqref{system-final}, \eqref{KdV system q}, respectively, for  the point vortex and  the surface variables.
The motion of the  point vortex is influenced by the free surface and the flat bottom.
In  the evolution of the surface variables, we observe the presence of the vorticity and of the vertical  derivative of the Green function evaluated on the free surface. It is interesting to note that even in the simplified case of an unperturbed KdV soliton travelling above the  vortex, its trajectories are bounded in time and their shape depends on the depth at which  the point vortex is located.

\section{Preliminaries}

\begin{figure}[h!]
\begin{center}
\fbox{\includegraphics[totalheight=0.3\textheight]{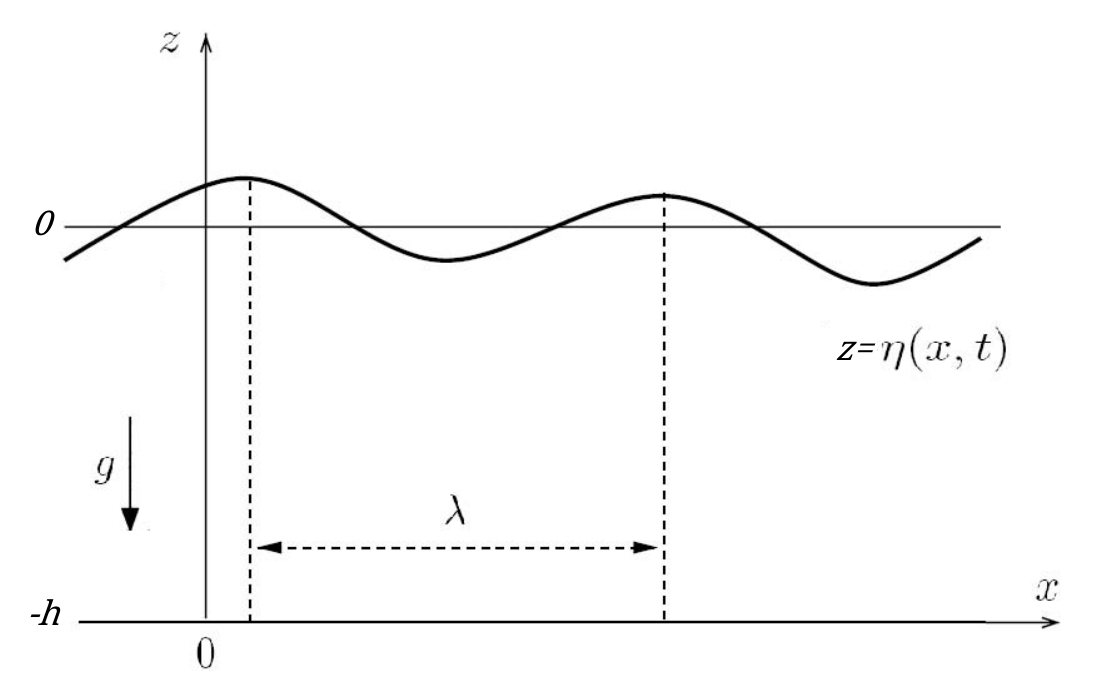}}
\caption{The fluid domain.}
\label{fig:setup}
\end{center}
\end{figure}

We consider a two-dimensional inviscid incompressible fluid in a constant
gravitational field. Let $x$ and $z$ be the horizontal and vertical coordinates, respectively. The fluid domain
$$\Omega_\eta\ =\ \{(x,z)\in\Bbb{R}^2:\  -h\ <\ z\ <\ \eta(x,t)\}$$
is bounded below by a flat rigid bottom
$$B\ =\ \{(x,z)\in\Bbb{R}^2:\ z\ =\ -h\}$$
and above by the free surface
$$S_\eta\ =\ \{(x,z)\in\Bbb{R}^2:\ z\ =\ \eta(x,t)\},$$
which we assume to be the graph of a function, see Fig. \ref{fig:setup}. Let $(u(x,z,t),v(x,z,t))$ be the velocity field.

The motion of the ideal fluid with a constant density $\rho$ set to 1, is described by Euler’s
equations
\begin{equation}
\begin{array}{cc}
u_t\ +\ uu_x\ +\ vu_z\ =\ -\ p_x& \\
 v_t\ +\ uv_x\ +\ vv_z\ =\ -\ p_z&-\ g\\
\end{array} \quad \textrm{in} \quad \Omega_{\eta},
 \label{Euler}
 \end{equation}
where  $p(x,z,t)$ denotes the pressure,
$g$ the constant gravitational acceleration in the negative $z$
direction, and the incompressibility condition
 \begin{align}
u_x\ +\ v_z\ =\ 0 \quad \textrm{in} \quad \Omega_{\eta}. \label{incomp}\end{align}
Complementing the equations of motion are the dynamic and the kinematic  boundary conditions
\begin{align}
p\ =\ p_{\textrm{atm}} \,  &\quad \textrm{ on }\quad  S_\eta
 \label{bc1}\\
  v\ =\ \eta_t\ +\ u\eta_x \, \, & \quad \textrm{ on }\quad S_\eta
 \label{bc2}\\
 v\ =\ 0 \, \,
& \quad  \textrm { on }\quad  B,
 \label{bc2'}
 \end{align}
with $p_{\textrm{atm}}$ the constant atmospheric pressure at the free surface.\\
We define now the scalar vorticity of the two-dimensional  flow by
\begin{align} \omega(x,z,t)\ =\ u_z\ -\  v_x\quad \textrm{ in } \quad \Omega_{\eta}.
\label{vorticity}\end{align}
Taking the curl of Euler equations \eqref{Euler}, we get the evolution of the vorticity in time
  \begin{align}
\omega_t\ +\ u\, \omega_x\ +\ v\, \omega_z\ =\ 0 \quad \textrm{ in } \quad \Omega_{\eta}.\label{omega_t}
\end{align}
Irrotational flow occurs when $\omega$  is zero everywhere, constant non-zero vorticity corresponds to a linear shear flow and non-constant vorticity
indicates highly sheared flows.
Throughout this paper we consider rotational flows of arbitrary nonconstant vorticity.

\section{Equations of motion reformulated with the aid of the Dirichlet-Neumann operator and the Green function in the free boundary domain}

For two-dimensional flows, the incompressibility condition \eqref{incomp} ensures the existence
of a stream function $\psi(x, z, t)$ - determined up to an additive term that depends
only on time -  by
\begin{align}
u\ =\ \psi_z, \quad v\ =\ -\psi_x.\label{psi}
\end{align}
Taking into account the kinematic boundary conditions \eqref{bc2}-\eqref{bc2'}, the stream functions
satisfies
\begin{equation}
 \frac{d}{dx} \psi (x, \eta(x,t),t)\ =\ -\eta_t \quad \textrm{ on }\quad S_\eta,
\end{equation}
and
\begin{align}
\psi_x(x,-h,t)\ =\ 0 \quad \textrm{ on }\quad B.
\end{align}
If the flow is irrotational, then the velocity can be expressed as the gradient of a scalar field, $\nabla \varphi$; one refers to $\varphi$ as the velocity potential. In our case, according to Weyl-Hodge theory (see \cite{Ebin&Marsden} for a summary and references), the velocity field
decomposes uniquely as
\begin{equation}\label{eq1}
  \begin{split}
     u\ =\ & \varphi_x\ +\ r_z\ \stackrel{\eqref{psi}}{=}\ \psi_z \\
      v\ =\ & \varphi_z\ -\ r_x\ \stackrel{\eqref{psi}}{=}\ -\psi_x
  \end{split}
\end{equation} where $\varphi(x,z,t)$ is the velocity potential, and $(r_z,-r_x)$ is a divergence free vector field, which is tangent to the free surface $S_{\eta}$.
The unit outward normal to the free surface is
\begin{align}{\bf n}\ =\ \frac{1}{\sqrt{1+\eta_x^2}}(-\eta_x, 1).\label{normal}
\end{align}
Thus, the function $r(x,z,t)$ has to satisfy the following condition
\begin{equation}\label{GF}
 r_x\ +\ r_z \eta_x\ =\ 0 \quad \implies \quad \frac{d}{dx}  r(x, \eta(x,t),t)\ =\ 0 \quad \textrm{ on }\quad S_\eta.
  \end{equation}
By \eqref{eq1},
\begin{align}
\omega\ =\ \ \Delta r\ =\ \ \Delta \psi \quad \textrm{ in } \quad \Omega_{\eta},\label{omega}
\end{align}
and the function $(\psi-r)$ is the harmonic conjugate of the function $\varphi$.\\
We will concentrate on waves with decay and so we assume that the functions  $r(x,z,t)$, $\eta(x,t)$, $\varphi(x,z,t)$ are fast-decaying functions as $|x| \to \infty$,  for all $t$, and on the boundaries we take
\begin{equation}
 r(x, \eta(x,t),t)\ =\ 0 \quad \textrm{ on }\quad S_\eta,\label{cond_r}
\end{equation}
 \begin{align}
 r(x,-h,t)\ =\ 0 \quad \textrm{ on }\quad B.\label{r_b}
 \end{align}
Let us now write the  Euler equations \eqref{Euler} in terms of the functions $\varphi$, $r$, and $\psi$.
With \eqref{eq1} and \eqref{omega} in view, we get
\begin{align}
\nabla\left(\varphi_t+\frac{1}{2}|\nabla\psi|^2+p+gz\right)& =\left(
\begin{array}{ll}
-r_{zt}\ -\ v\, \omega\\
\,\,\,\,  r_{xt}\ +\ u\, \omega
\end{array}\right)^T.\label{1}
\end{align}
From the vorticity equation \eqref{omega_t} and the incompressibility condition \eqref{incomp}, we have
\begin{align}r_{xt}\ +\ u\, \omega\ =\ -\partial_x^{-1}\left[(\ r_{zt}\ +\ v\, \omega)_z\right]\label{2}
\end{align}
and thus \eqref{1} becomes
\begin{align}
\nabla\left(\varphi_t\ +\ \frac{1}{2}|\nabla\psi|^2\ +\ p\ +\ gz\ +\  \partial_x^{-1}\left[\ r_{zt}\ +\ v\, \omega\right]\right)\ =\ 0.\label{3}
\end{align}
Hence,  at each instant $t$,
\begin{align*}
\varphi_t\ +\ \frac{1}{2}|\nabla\psi|^2\ +\  p\ +\ gz\ +\ \partial_x^{-1}\left[\ r_{zt}\ -\  \omega\, \psi_x\right]
\end{align*}
is constant throughout the fluid domain $\Omega_\eta$.
 By \eqref{bc1}, the pressure at the free surface is the constant atmospheric pressure, thus, by evaluating  at the free surface\footnote{We use subindex $S$ for evaluations at the free surface.} the expression above, we obtain that
$$
(\varphi_t)_S\ +\ \frac{1}{2}|\nabla\psi|_S^2\ +\ g\eta\ +\
\big[\partial_x^{-1}\left[\ r_{zt}\ -\  \omega\, \psi_x\right]\big]_S
$$
is constant on $S_\eta$. We absorb into the definition \eqref{eq1} of $\varphi$ a
suitable time-dependent term so that,  on $S_\eta$,
\begin{align}\label{bernoulli}
(\varphi_t)_S\ +\ \frac{1}{2}|\nabla\psi|_S^2\  +\ g\eta\ +\
\big[\partial_x^{-1}\left[\ r_{zt}\ -\  \omega\, \psi_x\right]\big]_S
=\ 0,
\end{align}
or in the form
\begin{align}
(\varphi_t)_S\ +\ \frac{1}{2}|\nabla\varphi|_S^2\ +\ \frac{1}{2}|\nabla r|_S^2\ +\ (\varphi_xr_{z} - r_x\varphi_z)_S
\qquad \qquad  &\nonumber\\
 +\ g\eta\ +\ \big[\partial_x^{-1}\left[\ r_{zt}\ -\  \omega\, \psi_x\right]\big]_S=\ 0.
\end{align}
Summing us, in terms of $\varphi(x,z,t)$, $r(x,z,t)$, $\eta(x,t)$, fast-decaying functions as $|x| \to \infty$,  for all $t$, we have the following equations
\begin{align}\label{eqs}
\hspace{-1.3cm}\begin{split}
\Delta \varphi\ =\ 0 \quad & \textrm{ in }\quad  \Omega_\eta \\
\Delta r\ =\  \omega \quad & \textrm{ in }\quad \Omega_\eta\\
(\varphi_t)_S\ +\  \frac{1}{2}|\nabla\psi|_S^2\
 +\ g\eta\ +\ \big[\partial_x^{-1}\left[\ r_{zt}\ -\  \omega\, \psi_x\right]\big]_S=\ 0 & \textrm{ on }\quad S_\eta\\
\eta_t\ =\ - \psi_x\Big|_S\ -\ \psi_z\Big|_S\eta_x\ =\ \varphi_z\Big|_S\ -\ \varphi_x\Big|_S\eta_x\ =\ \sqrt{1+\eta_x^2}\ \nabla\varphi \cdot \mathbf{n}  \quad & \textrm{ on }\quad S_\eta\\
 r_x\ +\ r_z\eta_x\ =\ 0 \quad & \textrm{ on }\quad S_\eta\\
 r\ =\ 0 \quad & \textrm{ on }\quad S_\eta\\
  \varphi_z\ =\ 0 \quad & \textrm{ on }\quad B\\
r\ =\ 0 \quad & \textrm{ on }\quad B.
\end{split}
\end{align}
Given initially an arbitrary vorticity $\omega$, with the evolution equation
\begin{align}\label{evol-omega}
\omega_t\ =\  \psi_x \omega_z\ -\ \psi_z\omega_x \quad \textrm{ in }\quad  \Omega_\eta,
\end{align}
 we will recast the problem \eqref{eqs}
in terms of the variables $\xi$,  $r$, $\eta$ and $\omega$, where  $\xi$ is
the
restriction of $\varphi$ to the free surface,
\begin{align}
\xi(x,t)\ :=\ \varphi(x,\eta(x,t),t) \quad \textrm{on} \quad S_{\eta}.\label{xi}
\end{align}
We introduce  the Dirichlet-Neumann operator $G$ associated to the domain $\Omega_\eta$, given by (see \cite{Craig})
\begin{align}\label{G}
G \xi :=\ \sqrt{1+\eta_x^2}\ (\nabla\varphi)_S \cdot \mathbf{n},
\end{align}
and hence we can write the forth equation in \eqref{eqs} as
\begin{align}\label{eta_t}
\eta_t\ =\  G\xi \quad \textrm{on} \quad S_{\eta}.
\end{align}
From \eqref{eqs}, the function $\varphi$ satisfies the boundary value problem
\begin{align}\label{system-varphi}
\begin{split}
\Delta \varphi\ =\  0 \quad & \textrm{ in }\quad \Omega_\eta\\
 \varphi_z\ =\ 0 \quad & \textrm{ on }\quad B\\
\varphi\ =\ \xi \quad & \textrm{ on }\quad S_\eta.
\end{split}\end{align}
%Under reasonable regularity assumptions on $\eta$ and $\xi$, specified by rigorous mathematical analysis  (see, for example, the Chapters 1 and 2 in \cite{lannes}), the problem above
%has a unique solution. More precisely, cf. Corollary 2.44 in \cite{lannes}, for any  $s_0>\frac{1}{2}$,  $\eta\in H^{s_0+1}(\mathbb{R})$ and $\xi\in \dot{H}^{\frac{3}{2}}(\mathbb{R})$, there exists a unique solution $\varphi\in H^2(\Omega_\eta)$ to the boundary value problem \eqref{system-varphi}.\\
We take a basis of elementary harmonic functions in the strip $-h<z<H$, where $H$=const is the upper bound of $\eta$, satisfying the first  boundary  condition in \eqref{system-varphi}, that is
\begin{align}
\varphi_k(x,z)=\mathfrak{A}(k)\, \cosh\big((z+h)k\big){\rm e}^{ikx}.\
\end{align}
Then, we consider a subclass of solutions to \eqref{system-varphi}, regular in the strip $-h<z<H$, given by
\begin{align}\label{31}
\varphi(x, z)\ &= \  \frac{1}{\pi}\int_{-\infty}^\infty  \Big[
\mathfrak{A}(k)\, \cosh\big((z+h)k\big)\Big]\,  {\rm e}^{ikx}\,  dk,
\end{align}
for a suitable choice of $\mathfrak{A}(k)$, for example,
\begin{align}\label{cond1}
|\mathfrak{A}(k)|\leq C{\rm e}^{-a|k|}, \quad a>h+H>0,\, \,  C>0\, \, \textrm{const}.
\end{align}
The  existence and uniqueness of the solutions to the boundary value problem \eqref{system-varphi}  is studied, for example, in \cite{lannes}.\\
 We introduce the operator
\begin{align}\label{D}
\begin{split}
D=- i\, \partial _x,\quad \partial_x= i\, D,\quad  \partial_x^{-1}= -i\, D, \quad D^{-1}=i\, \partial_x^{-1},
\end{split}
\end{align}
 and taking into account that for an arbitrary function $f(x)$ the operator ${\rm e}^{\alpha D}$, with $\alpha\in\Bbb{R}$, applied to $f$ is
${\rm e}^{\alpha D} f(x)\ =\ {\rm e}^{-i\alpha\partial_x}f(x)\ =\ f(x-i\alpha)$, we have
\begin{align}\label{Dexp}
\begin{split}
{\rm e}^{-(z+h)D}{\rm e}^{ikx}\ &=\ {\rm e}^{ik[x+i(z+h)]}\ =\  {\rm e}^{ikx}{\rm e}^{-(z+h)k}.
\end{split}
\end{align}
Hence, the function $\varphi$ in \eqref{31} can be written as
\begin{align}\label{32}
\varphi(x, z)\ &= \   \cosh\big((z+h)D\big) \int_{-\infty}^\infty
\frac{1}{\pi}\mathfrak{A}(k)\,  {\rm e}^{ikx}\,  dk,
\end{align}
where the operator
 \begin{align}
\cosh\big((z+h)D\big)\ =\  1\ +\ \frac{(z+h)^2D^2}{2!}\ + \ \frac{(z+h)^4D^4}{4!}\ +\ ...
\end{align}
By the second boundary conditions in \eqref{system-varphi}, we have
\begin{align}
\xi(x)\ =\ : \cosh\big((\eta+h)D\big): \int_{-\infty}^\infty
\frac{1}{\pi}\mathfrak{A}(k)\,  {\rm e}^{ikx}\,  dk,\nonumber
\end{align}
 thus, we get
\begin{align}\label{varphi}
\varphi(x,z)\ =\ \cosh\big((z+h)D\big)\big[:\cosh\big((\eta+h)D\big):\big]^{-1}\ \xi,
\end{align}
where  $:\cosh\big((\eta+h)D\big):$ is the normal ordering of the operator\\ $\cosh\big((\eta+h)D\big)$,  that is, in its expression
\begin{align}
:\cosh\big((\eta+h)D\big):\ =\  1\ +\ \frac{(\eta+h)^2D^2}{2!}\ + \ \frac{(\eta+h)^4D^4}{4!}\ +\ ... ,
\end{align}
first are the powers of the function $\eta(x,t)+h$  which depends on $x$ and then the powers of the derivative operator  $D$.

 \noindent Let us focus now on the other boundary value problem
\begin{align}\label{system-r}
\begin{split}
\Delta r\ =\  \omega \quad & \textrm{ in }\quad \Omega_\eta\\
 r\ =\ 0 \quad & \textrm{ on }\quad B\\
 r\ =\ 0 \quad & \textrm{ on }\quad S_\eta.
\end{split}
\end{align}
We seek to solve this problem by using a Green function, that is, the solution  $\Gamma(\mathbf{x}, \mathbf{x}_0)$ to  the problem
\begin{align}\label{Gamma_omega}
\begin{split}
\Delta \Gamma\ &=\ \delta(\mathbf{x}-\mathbf{x}_0)  \quad  \textrm{ in }\quad \Omega_\eta \\
\Gamma\ &=\ 0 \quad \quad\qquad\,\,\,  \textrm{ on } \quad  B\cup S_\eta,
\end{split}
\end{align}
where $\delta$ is the Dirac delta function, $\mathbf{x}=(x,z)\in\Omega_\eta$ and $\mathbf{x}_0=(x_0,z_0)$ is a fixed point in $\Omega_\eta$.
Then, by  Green's  representation formula (see, for example,  \cite{Haberman}), the solution to the problem \eqref{system-r} can be written as
\begin{align}\label{superposition}
r(\mathbf{x})\ =\  \Delta^{-1} \omega\ :=\  \int_{\Omega_\eta} \Gamma(\mathbf{x}, \mathbf{x_0})\ \omega(\mathbf{x_0})\ d\mathbf{x}_0.
\end{align}
In order to construct the Green function $\Gamma(\mathbf{x},\mathbf{x}_0)$ for our domain $\Omega_\eta\subset \Bbb{R}^2$, we use the fundamental Green function  $\Gamma_0(\mathbf{x},\mathbf{x}_0)$ in the whole plane, that is, the function which satisfies the equation
$$\Delta \Gamma_0\ =\ \delta (\mathbf{x}-\mathbf{x}_0)  \quad  \textrm{ in } \quad  \Bbb{R}^2,$$
for any point $\mathbf{x}_0\in \Bbb{R}^2$.  This function  has the expression  (see, for example, \cite{Haberman})
\begin{align}\label{Gamma_0}
\Gamma_{0}(\mathbf{x},\mathbf{x}_0)\ =\  \frac{1}{4\pi}\ln \Big((x-x_0)^2\ +\ (z-z_0)^2\Big).
\end{align}
We set
\begin{align}\label{mathcalF}
\Gamma(\mathbf{x},\mathbf{x}_0)\ =\  \Gamma_0(\mathbf{x},\mathbf{x}_0)\ -\  \mathfrak{F}(\mathbf{x},\mathbf{x}_0),
\end{align}
where  $\mathfrak{F}(\mathbf{x},\mathbf{x}_0)$  is a harmonic function in $\Omega_\eta$  so that the Green function $\Gamma$ satisfies the boundary condition in \eqref{Gamma_omega}, that is,
\begin{align}\label{F}
\begin{split}
\Delta \mathfrak{F}\ &=\ 0\quad\quad\quad\quad  \qquad\qquad\quad \textrm{ in }\quad \Omega_\eta \\
\mathfrak{F}\ &=\   \Gamma_0(x,-h,\mathbf{x}_0)\ \quad\quad\quad \,\, \textrm{ on }\quad B\\
\mathfrak{F}\ &=\    \Gamma_0(x,\eta(x,t),\mathbf{x}_0)  \quad\quad\,\,  \textrm{ on }\quad S_\eta.
\end{split}
\end{align}
Similar to the assumptions made in relation to the boundary value problem \eqref{system-varphi}, we now assume the
existence of a subclass of solutions $\mathfrak{F}(\mathbf{x}, \mathbf{x}_0)$, regular in the strip $-h<z<H$.\\
With the help of the operator $D$ defined in \eqref{D}, \eqref{Dexp},
we write the general solution to the above  equation  into the form
\begin{align}
\mathfrak{F}(\mathbf{x}, \mathbf{x}_0)\ =&\   {\rm e}^{-(z+h)D}\int_{-\infty}^\infty
\frac{\mathcal{A}(k,\mathbf{x}_0)}{2\pi}  \,  {\rm e}^{ikx}\,  dk\ +\  {\rm e}^{(z+h)D}\int_{-\infty}^\infty
\frac{\mathcal{B}(k,\mathbf{x}_0)}{2\pi} \,  {\rm e}^{ikx}\,  dk.\nonumber
\end{align}
From the   boundary conditions in \eqref{F} we get
\begin{align}
\int_{-\infty}^\infty
\frac{\mathcal{A}(k,\mathbf{x}_0)}{2\pi}  \,  {\rm e}^{ikx}\,  dk\ =&\   - \frac{1}{2}\left[:\sinh\big((\eta+h)D\big):\right]^{-1}
\Gamma_0(x,\eta,\mathbf{x}_0)\nonumber\\
& +\  \frac{1}{2}\left[: \sinh\big((\eta+h)D\big):\right]^{-1}  :{\rm e}^{(\eta+h)D}:
  \Gamma_0(x,-h,\mathbf{x}_0)
\nonumber
\end{align}
\begin{align}
\int_{-\infty}^\infty
\frac{\mathcal{B}(k,\mathbf{x}_0)}{2\pi} \,  {\rm e}^{ikx}\,  dk\ =&\  \frac{1}{2}  \left[: \sinh\big((\eta+h)D\big):\right]^{-1}
\Gamma_0(x,\eta,\mathbf{x}_0)\nonumber\\
& -\ \frac{1}{2} \left[:\sinh\big((\eta+h)D\big):\right]^{-1} : {\rm e}^{-(\eta+h)D}:
  \Gamma_0(x,-h,\mathbf{x}_0).
\nonumber
\end{align}
Finally, we conclude that
\begin{align}
 \mathfrak{F}(\mathbf{x}, \mathbf{x}_0)\ =&\  \sinh\big((z+h)D\big)\left[: \sinh\big((\eta+h)D\big):\right]^{-1} \Gamma_0(x,\eta,\mathbf{x}_0) \nonumber\\
&+\  \frac{1}{2}{\rm e}^{-(z+h)D} \left[:\sinh\big((\eta+h)D\big):\right]^{-1}:{\rm e}^{(\eta+h)D}:  \Gamma_0(x,-h,\mathbf{x}_0)\nonumber\\
&-\  \frac{1}{2} {\rm e}^{(z+h)D} \left[:\sinh\big((\eta+h)D\big):\right]^{-1}:{\rm e}^{-(\eta+h)D}:  \Gamma_0(x,-h,\mathbf{x}_0) \nonumber\\
=&\  \sinh\big((z+h)D\big)\left[: \sinh\big((\eta+h)D\big):\right]^{-1} \Gamma_0(x,\eta,\mathbf{x}_0) \nonumber\\
&-\  :\sinh\big((z-\eta)D\big): \left[:\sinh\big((\eta+h)D\big):\right]^{-1}  \Gamma_0(x,-h,\mathbf{x}_0).\nonumber
\end{align}
Therefore, the Green function \eqref{mathcalF} has the expression
\begin{align}\label{Gamma}
 \Gamma(\mathbf{x},\mathbf{x}_0)\
 =& \   \Gamma_0(\mathbf{x},\mathbf{x}_0)\ -\
   \sinh\big((z+h)D\big)\left[: \sinh\big((\eta+h)D\big):\right]^{-1} \Gamma_0(x,\eta,\mathbf{x}_0) \nonumber\\
& +\ :\sinh\big((z-\eta)D\big):\left[:\sinh\big((\eta+h)D\big):\right]^{-1} \Gamma_0(x,-h,\mathbf{x}_0),
\end{align}
and
\begin{align}\label{r}
r(\mathbf{x})\ =&\   \int_{\Omega_\eta}  \Gamma(\mathbf{x},\mathbf{x}_0)\ \omega(\mathbf{x_0})\ d\mathbf{x}_0.
\end{align}
From \eqref{eq1}, \eqref{D}, \eqref{varphi}, {\eqref{r} and the relation
\begin{align}
D^{-1} \sinh\big((z+h)D\big)\ =\  \sinh\big((z+h)D\big) D^{-1},\nonumber
\end{align}
 we obtain that
\begin{align}\label{50}
\psi(x,z)\
=\  &  \int_{\Omega_\eta}  \Gamma(\mathbf{x},\mathbf{x}_0)\ \omega(\mathbf{x_0})\ d\mathbf{x}_0\nonumber\\&+i\   \sinh\big((z+h)D\big) \, \left[:\cosh\big((\eta+h)D\big):\right]^{-1}\xi.
\end{align}
Therefore,
\begin{align}
\psi_x\ =&\  \int_{\Omega_\eta}  \Gamma_x(\mathbf{x},\mathbf{x}_0)\ \omega(\mathbf{x_0})\ d\mathbf{x}_0
\nonumber\\
&-   \sinh\big((z+h)D\big) D \left[:\cosh\big((\eta+h)D\big):\right]^{-1}\xi\end{align}
\begin{align}
\psi_x\Big|_S\ =&\  \Big(\int_{\Omega_\eta}  \Gamma_x(\mathbf{x},\mathbf{x}_0)\ \omega(\mathbf{x_0})\ d\mathbf{x}_0\Big)_S
\nonumber\\
&-   :\sinh\big((\eta+h)D\big): D \left[:\cosh\big((\eta+h)D\big):\right]^{-1}\xi\end{align}
\begin{align}
\psi_z\ =&\  \int_{\Omega_\eta}  \Gamma_z(\mathbf{x},\mathbf{x}_0)\ \omega(\mathbf{x_0})\ d\mathbf{x}_0
\nonumber\\
&+i\   \cosh\big((z+h)D\big)D\left[:\cosh\big((\eta+h)D\big):\right]^{-1}\xi\end{align}
\begin{align}
\psi_z\Big|_s\ =&\  \Big(\int_{\Omega_\eta}  \Gamma_z(\mathbf{x},\mathbf{x}_0)\ \omega(\mathbf{x_0})\ d\mathbf{x}_0\Big)_S
\nonumber\\
&\hspace{2cm}+i\  : \cosh\big((\eta+h)D\big): D\left[:\cosh\big((\eta+h)D\big):\right]^{-1}\xi\nonumber\\
=& \Big(\int_{\Omega_\eta}  \Gamma_z(\mathbf{x},\mathbf{x}_0)\ \omega(\mathbf{x_0})\ d\mathbf{x}_0\Big)_S\ +\  \xi_x
\nonumber\\
&\hspace{0.5cm}\ -\, \eta_x:\sinh\big((\eta+h)D\big): D\left[:\cosh\big((\eta+h)D\big):\right]^{-1} \xi
\end{align}
In view of the forth equation in  \eqref{eqs},  from  \eqref{eta_t}, we have
\begin{align}\label{verification}
\psi_x\Big|_S+\eta_x\psi_z\Big|_s=-G\xi,
\end{align}
that is,
\begin{align}\label{51}
&  \Big(\int_{\Omega_\eta}  \Gamma_x(\mathbf{x},\mathbf{x}_0)\ \omega(\mathbf{x_0})\ d\mathbf{x}_0\Big)_S+\eta_x\Big(\int_{\Omega_\eta}  \Gamma_z(\mathbf{x},\mathbf{x}_0)\ \omega(\mathbf{x_0})\ d\mathbf{x}_0\Big)_S
\nonumber\\
&+\eta_x\xi_x\ -   :\sinh\big((\eta+h)D\big): D \left[:\cosh\big((\eta+h)D\big):\right]^{-1}\xi\nonumber\\
&-\ \eta_x^2:\sinh\big((\eta+h)D\big): D\left[:\cosh\big((\eta+h)D\big):\right]^{-1} \xi\nonumber\\
&=\eta_x\xi_x-(1+\eta_x^2) : \sinh\big((\eta+h)D\big): D\, \left[: \cosh\big((\eta+h)D\big):\right]^{-1} \xi,
\end{align}
where, by the fifth equation in \eqref{eqs}, the sum of the above integrals is zero.

\noindent We underline the fact that the Green function in \eqref{Gamma} depends explicitly on the free surface $\eta$.
From \eqref{r}, the evolution equation in time for the function $r$ is
\begin{align}\label{r_t}
r_t(\mathbf{x})\ =&\ \mathcal{R}(x,z),
\end{align}
where
\begin{align}
\mathcal{R}(x,z):=& \int_{-\infty}^{\infty} \int_{-h}^{\eta(x_0,t)}\Big[ \Gamma_t(\mathbf{x},\mathbf{x}_0)\ \omega(\mathbf{x_0})\ +\   \Gamma(\mathbf{x},\mathbf{x}_0)\ \omega_t(\mathbf{x_0}) \Big]dz_0\ dx_0
 \nonumber\\
&\ +\int_{-\infty}^\infty \Gamma(\mathbf{x},x_0,\eta(x_0,t))\  \omega(x_0,\eta(x_0,t))\  \eta_t(x_0,t)\  dx_0\nonumber\\
\stackrel{\eqref{evol-omega}, \eqref{eta_t}}{=}&\ \int_{\Omega_\eta} \Big[   \Gamma_t(\mathbf{x},\mathbf{x}_0)\ \omega(\mathbf{x_0})\ + \Gamma(\mathbf{x},\mathbf{x}_0)\  ( \psi_{x_0} \omega_{z_0}\ -\ \psi_{z_0}\omega_{x_0} ) \Big]d\mathbf{x}_0\nonumber\\
&\ +\int_{-\infty}^\infty \Gamma(\mathbf{x},x_0,\eta(x_0,t))\  \omega(x_0,\eta(x_0,t))\  G \xi(x_0) \  dx_0.
\end{align}
From \eqref{Gamma}, and by $(\Gamma_0)_t(\mathbf{x},\mathbf{x}_0)=0$, for $\mathbf{x}_0$  fixed in the fluid domain,   the time derivative of $\Gamma$, which depends explicitly on $\eta$, is
\begin{align}
 \Gamma_t(\mathbf{x},\mathbf{x}_0)\ =\   \Gamma_\eta(\mathbf{x},\mathbf{x}_0)\, \eta_t\ \stackrel{\eqref{eta_t}}{=}   \Gamma_\eta(\mathbf{x},\mathbf{x}_0) \, G \xi(x) .
\end{align}
It remains to look at the evolution equation for $\xi$ and at the   Bernoulli relation on the free surface, that is, the third relation in the system \eqref{eqs}.
First, by differentiating \eqref{xi} with respect to $x$ we get
\begin{align}
\xi_x\ =\  (\varphi_x)_S\ + \ (\varphi_z)_S\eta_x,
\end{align}
and  from \eqref{G}, we have the system
\begin{align}\label{systemG}
\begin{split}
 (\varphi_x)&_S\ + \ (\varphi_z)_S\eta_x\ =\ \xi_x\\
 - (\varphi_x)&_S\eta_x\  +\  (\varphi_z)_S =G\xi.
\end{split}
\end{align}
The system \eqref{systemG} can be solved as
\begin{align}\label{systemGsol}
\begin{split}
(\varphi_x)_S\  & =\ \frac{1}{1+\eta_x^2}\ \Big( \xi_x - \eta_x G \xi \Big) \\
(\varphi_z)_S\   &=\  \frac{1}{1+\eta_x^2}\ \Big( G\xi  + \eta_x \xi_x\Big).
\end{split}
\end{align}
On the other hand, taking into account \eqref{varphi},
\begin{align}
(\varphi_z)_S\ =\  :\sinh\big((\eta+h)D\big): D\, \left[: \cosh\big((\eta+h)D\big):\right]^{-1} \xi.
\end{align}
Therefore, by \eqref{systemGsol}, we obtain the following formula for the Dirichlet-Neumann operator in terms of the free surface
and the flat bottom
\begin{align}\label{G1}
G \xi\ =&\  -\eta_x\xi_x\ +\ (1+\eta_x^2)\, : \sinh\big((\eta+h)D\big): D\, \left[: \cosh\big((\eta+h)D\big):\right]^{-1} \xi,
\end{align}
formula that agrees with \eqref{51}.
In Section 5, we will compare
the above  new explicit  expression for the Dirichlet-Neumann operator with the commonly used one as a Taylor series   \cite{Craig, Craig&Groves, CS}, see the formula  \eqref{G3}.\\
Differentiation \eqref{xi} with respect to $t$ we get
\begin{align}
\xi_t\ =\  (\varphi_t)_S\ + \ (\varphi_z)_S\eta_t,\nonumber
\end{align}
and in view of  the third relation in the system \eqref{eqs},
\begin{align}
\xi_t\ =& \ (\varphi_z)_S\ \eta_t\  -\frac{1}{2}|\nabla\psi|_S^2\  -g\eta\ -\   \big[\partial_x^{-1}\left[\ r_{zt}\ -\  \omega\, \psi_x\right]\big]_S.\end{align}
According to \eqref{eta_t}, \eqref{r_t} and \eqref{systemGsol}, we obtain
\begin{align}
\xi_t\
=&\ \frac{\eta_x\,  \xi_x}{1+\eta_x^2}\,  G\xi  +\  \frac{1}{1+\eta_x^2} \,  \big(G\xi\big)\cdot \big(G\xi\big)\  -\frac{1}{2}|\nabla\psi|_S^2\  -g\eta \nonumber\\
&\  -\   \big[\partial_x^{-1}\left[\ \mathcal{R}_z(x,z)\ -\  \omega\, \psi_x\right]\big]_S.
\end{align}
For the last term above we can  calculate
\begin{align}
\frac{d}{dx}\Big[\big[\partial_x^{-1}\left[\ r_{zt}\ -\  \omega\, \psi_x\right]\big]_S\Big]\
 =& \ (r_{zt})_S\ -\ (\omega\psi_x)_S\ +\
\big[\partial_x^{-1}\left[(\ r_{zt}\ -\  \omega\, \psi_x)_z\right]\big]_S
\eta_x\nonumber\\
\stackrel{\eqref{2}, \eqref{eqs}_4}{=}& (r_{zt})_S\ -\ (r_{xt})_S\eta_x\ +\ \omega_S\eta_t.
\end{align}

\noindent Summing up, we have proved that:

\begin{theorem}
For an arbitrary vorticity $\omega(x,z,t)$,  the nonlinear governing equations for the  wave-vorticity interactions, can be written  in terms of the variables   $\omega(x,z,t)$, $r(x,z,t)$,   $\eta(x,t)$   and $\xi(x,t)$, in the following way:
\begin{align}\label{eqs1}
\begin{split}
\omega_t\ =&\  \psi_x \omega_z\ -\ \psi_z\omega_x \qquad \qquad\qquad \hspace{4.3cm} \textrm{ in }\,\,  \Omega_\eta\\
r_t\ =&\ \mathcal{R}(x,z) \hspace{7.8cm}  \textrm{ in }\,\,  \Omega_\eta\\
\eta_t\ =& \ G\xi \hspace{8.7cm} \textrm{on} \,\, S_{\eta}\\
\xi_t\ =&\   \frac{\eta_x\,  \xi_x}{1+\eta_x^2}\,  G\xi  +\  \frac{1}{1+\eta_x^2} \,  \big(G\xi\big)\cdot \big(G\xi\big)\  -\frac{1}{2}|\nabla\psi|_S^2\  -g\eta\\
&\ +\  \partial_x^{-1}\Big[\eta_x \mathcal{R}_x(x,z)\Big|_S\ -\    \mathcal{R}_z(x,z)\Big|_S \ -\  \omega_S(G\xi)\Big]\hspace{1.2cm} \textrm{on} \,\, S_{\eta}
\end{split}
\end{align}
 where
\begin{align}\label{psi2}
\psi(x,z)\ =&\  \int_{\Omega_\eta}  \Gamma(\mathbf{x},\mathbf{x}_0)\ \omega(\mathbf{x_0})\ d\mathbf{x}_0
\nonumber\\
&+i\   \sinh\big((z+h)D\big) \, \left[:\cosh\big((\eta+h)D\big):\right]^{-1}\xi\end{align}
\begin{align}\label{61}
 \Gamma(\mathbf{x},\mathbf{x}_0)\ =&\    \Gamma_0(\mathbf{x},\mathbf{x}_0) \nonumber\\
&-\
   \sinh\big((z+h)D\big)\left[: \sinh\big((\eta+h)D\big):\right]^{-1} \Gamma_0(x,\eta,\mathbf{x}_0) \nonumber\\
&+\ :\sinh\big((z-\eta)D\big):\left[:\sinh\big((\eta+h)D\big):\right]^{-1} \Gamma_0(x,-h,\mathbf{x}_0)\end{align}
\begin{align}\label{62}
 \Gamma_0(\mathbf{x},\mathbf{x}_0)\ =\ \frac{1}{4\pi}\ln \Big((x-x_0)^2\ +\ (z-z_0)^2\Big)
\end{align}
\begin{align}\label{psi_x}
\psi_x\ =&\  \int_{\Omega_\eta}  \Gamma_x(\mathbf{x},\mathbf{x}_0)\ \omega(\mathbf{x_0})\ d\mathbf{x}_0
\nonumber\\
&-   \sinh\big((z+h)D\big) D \left[:\cosh\big((\eta+h)D\big):\right]^{-1}\xi\end{align}
\begin{align}\label{psi_x_S}
\psi_x\Big|_S\ =&\  \Big(\int_{\Omega_\eta}  \Gamma_x(\mathbf{x},\mathbf{x}_0)\ \omega(\mathbf{x_0})\ d\mathbf{x}_0\Big)_S
\nonumber\\
&-   :\sinh\big((\eta+h)D\big): D \left[:\cosh\big((\eta+h)D\big):\right]^{-1}\xi\end{align}
\begin{align}\label{psi_z}
\psi_z\ =&\  \int_{\Omega_\eta}  \Gamma_z(\mathbf{x},\mathbf{x}_0)\ \omega(\mathbf{x_0})\ d\mathbf{x}_0
\nonumber\\
&+i\   \cosh\big((z+h)D\big)D\left[:\cosh\big((\eta+h)D\big):\right]^{-1}\xi\end{align}

\begin{align}\label{psi_z_S}
\psi_z\Big|_s\ =& \Big(\int_{\Omega_\eta}  \Gamma_z(\mathbf{x},\mathbf{x}_0)\ \omega(\mathbf{x_0})\ d\mathbf{x}_0\Big)_S\ +\  \xi_x
\nonumber\\
&\hspace{0.5cm}\ -\, \eta_x:\sinh\big((\eta+h)D\big): D\left[:\cosh\big((\eta+h)D\big):\right]^{-1} \xi
\end{align}
\begin{align}\label{R(x,z)}
\mathcal{R}(x,z)
=&\ \int_{\Omega_\eta} \Big[ (G\xi)\, \Gamma_\eta(\mathbf{x},\mathbf{x}_0)\ \omega(\mathbf{x_0})\ +  \Gamma(\mathbf{x},\mathbf{x}_0)\  ( \psi_{x_0} \omega_{z_0}\ -\ \psi_{z_0}\omega_{x_0} ) \Big]d\mathbf{x}_0\nonumber\\
&\ +\int_{-\infty}^\infty \Gamma(\mathbf{x},x_0,\eta(x_0))\  \omega(x_0,\eta(x_0))\  G\xi (x_0) \  dx_0
\end{align}
\begin{align}\label{G2}
G \xi\ =&\  -\eta_x\xi_x\ +\ (1+\eta_x^2): \sinh\big((\eta+h)D\big): D\, \left[: \cosh\big((\eta+h)D\big):\right]^{-1} \xi
\end{align}
\end{theorem}

\section {Point vortex}

\noindent In this section we suppose that the vorticity
in the problem comes from a  point
vortex submerged beneath the free surface. Let us consider that
\begin{align}\label{omega*}
\omega\ =\ \omega^*\, \delta(\mathbf{x}-\mathbf{q})\ =\ \omega^*\delta(x-q_1)\delta(z-q_2),
\end{align}
where $\omega^*$ is a constant and $\mathbf{q}=(q_1(t),q_2(t))\in \Omega_\eta$ the position of the point vortex. Then,
\begin{align}
\omega_t\ =\ \omega^*\Big(
-\delta'(x-q_1)\delta(z-q_2)\dot{q_1}\ -\  \delta(x-q_1)\delta'(z-q_2)\dot{q_2}
\Big)
\end{align}
and the first equation in
 \eqref{eqs1} becomes
\begin{align}\label{78}
\begin{split}
\dot{q_1}\ =&\ \psi_z\Big|_{\mathbf{x}=\mathbf{q}}\\
 \dot{q_2}\ =& -\psi_x\Big|_{\mathbf{x}=\mathbf{q}}.
\end{split}
\end{align}
By  \eqref{psi2},  \eqref{psi_x} and \eqref{psi_z} we get
\begin{align}
\psi(x,z)\ =&\  \omega^*  \Gamma(\mathbf{x},\mathbf{q})+i\   \sinh\big((z+h)D\big) \, \left[:\cosh\big((\eta+h)D\big):\right]^{-1}\xi
\end{align}
\begin{align}
\psi_x\ =&\  \omega^*  \Gamma_x(\mathbf{x},\mathbf{q})
-   \sinh\big((z+h)D\big) D \left[:\cosh\big((\eta+h)D\big):\right]^{-1}\xi
\end{align}
\begin{align}
\psi_z\ =&\ \omega^*  \Gamma_z(\mathbf{x},\mathbf{q})+i\   \cosh\big((z+h)D\big)D\left[:\cosh\big((\eta+h)D\big):\right]^{-1}\xi .\end{align}
Hence, the equations \eqref{78} become
\begin{align}\label{eqsq}
\begin{split}
\dot{q_1}\
=&\  \omega^*\Gamma_z(\mathbf{x},\mathbf{q})\Big|_{\mathbf{x}=\mathbf{q}}\\
&+i\  \Big[ \cosh\big((z+h)D\big)D\left[:\cosh\big((\eta+h)D\big):\right]^{-1}\xi\Big]_{\mathbf{x}=\mathbf{q}}\\
 \dot{q_2}\
 =& -\omega^* \Gamma_x(\mathbf{x},\mathbf{q})\Big|_{\mathbf{x}=\mathbf{q}}\\
&+  \Big[\sinh\big((z+h)D\big) D\left[:\cosh\big((\eta+h)D\big):\right]^{-1}\xi\Big]_{\mathbf{x}=\mathbf{q}}.
\end{split}
\end{align}
Close to the vortex, that is,  $\mathbf{x}\rightarrow\mathbf{q}$, the Green function  \eqref{61}-\eqref{62} diverges as a logarithm.
 We  accept, as in \cite{M&Pulvirenti, Meleshko}, that there is no self-interaction.
By assuming that a single vortex does not move under the action of its own field,
in the system \eqref{eqsq} we  will  exclude  the singular terms.

\noindent The third equation in \eqref{eqs1} is not affected by the vorticity.\\
 Let us look what becomes the expression of the function $\mathcal{R}$ from \eqref{R(x,z)}
\begin{align}
\mathcal{R}(x,z)
=&\  \omega^*\eta_t\Gamma_\eta(\mathbf{x}, \mathbf{q})\ \nonumber\\
& +\omega^*\int_{\Omega_\eta} \Gamma(\mathbf{x},\mathbf{x}_0)\Big[  \psi_{x_0}
 \delta(x_0-q_1)\delta'(z_0-q_2)
 - \psi_{z_0}\delta'(x_0-q_1)\delta(z_0-q_2)  \Big]d\mathbf{x}_0\nonumber\\
&\ +\omega^*\int_{-\infty}^\infty \Gamma(\mathbf{x},x_0,\eta(x_0))\  \delta(x_0-q_1)\delta(\eta(x_0)-q_2)  G \xi(x_0) \  dx_0.\nonumber
\end{align}
Taking into account the fact that the point vortex can not be on the free surface, the last term above is zero, then,
\begin{align}
\mathcal{R}(x,z)\ =&\  \omega^*\eta_t\Gamma_\eta(\mathbf{x}, \mathbf{q})+\omega^*\Big(-\Gamma_{q_2}(\mathbf{x},\mathbf{q}) \psi_{q_1}(\mathbf{q}) -
\Gamma(\mathbf{x},\mathbf{q}) \psi_{q_1q_2}(\mathbf{q})\nonumber\\
&\hspace{3.5cm}\,\,+ \Gamma_{q_1}(\mathbf{x},\mathbf{q})\psi_{q_2}(\mathbf{q})+\Gamma(\mathbf{x},\mathbf{q})\psi_{q_1q_2}(\mathbf{q})\Big) \nonumber\\
=&\  \omega^*\eta_t\Gamma_\eta(\mathbf{x}, \mathbf{q})+\omega^*\Big(-\Gamma_{q_2}(\mathbf{x},\mathbf{q}) \psi_{q_1}(\mathbf{q})  + \Gamma_{q_1}(\mathbf{x},\mathbf{q})\psi_{q_2}(\mathbf{q})\Big),\nonumber\\
\stackrel{\eqref{78}}{=}&\  \omega^*\eta_t\Gamma_\eta(\mathbf{x}, \mathbf{q})+\omega^*\Big(\Gamma_{q_2}(\mathbf{x},\mathbf{q})\dot{q}_2  + \Gamma_{q_1}(\mathbf{x},\mathbf{q})\dot{q}_1\Big).
\end{align}
By the third equation in \eqref{eqs1},
the second equation in \eqref{eqs1} takes the form
\begin{align}
r_t\ =\  \omega^*\Gamma_\eta(\mathbf{x}, \mathbf{q})\, G\xi+\omega^*\Big(\Gamma_{q_1}(\mathbf{x},\mathbf{q})\dot{q}_1  + \Gamma_{q_2}(\mathbf{x},\mathbf{q})\dot{q}_2\Big)  \quad \textrm{ in } \,\, \Omega_\eta.
\end{align}
 From \eqref{omega*}, and the fact that  we do not allow point vortices on the free surface, we get $\omega_S=0$. Thus,  the forth equation in \eqref{eqs1} becomes
\begin{align}
\xi_t\
 =&
\  \frac{\eta_x\,  \xi_x}{1+\eta_x^2}\,  G\xi  +\  \frac{1}{1+\eta_x^2} \,  \big(G\xi\big)\cdot \big(G\xi\big)\  -\frac{1}{2}|\nabla\psi|_S^2\  -g\eta\nonumber\\
&\ +\  \partial_x^{-1}\Big[\eta_x \mathcal{R}_x(x,z)\Big|_S\ -\    \mathcal{R}_z(x,z)\Big|_S \Big]\hspace{1.2cm} \textrm{on} \,\, S_{\eta}.
\end{align}
\noindent Concluding,  \textit{the system for the surface variables and the point vortex is }
\begin{align}\label{eqs_point_v}
 \begin{split}
\dot{q_1}\
=&\  i  \Big[ \cosh\big((z+h)D\big)D\left[:\cosh\big((\eta+h)D\big):\right]^{-1}\xi\Big]_{\mathbf{x}=\mathbf{q}}\\
 \dot{q_2}\
 =&\,\,\,\,  \Big[\sinh\big((z+h)D\big) D\left[:\cosh\big((\eta+h)D\big):\right]^{-1}\xi\Big]_{\mathbf{x}=\mathbf{q}}\\
r_t\ =&\  \mathcal{R}(x,z)\\
\eta_t\ =& \ G\xi \\
\xi_t\ =&\   \frac{\eta_x\,  \xi_x}{1+\eta_x^2}\,  G\xi  +\  \frac{1}{1+\eta_x^2} \,  \big(G\xi\big)\cdot \big(G\xi\big)\  -\frac{1}{2}|\nabla\psi|_S^2\  -g\eta\\
&\ +\  \partial_x^{-1}\Big[\eta_x \mathcal{R}_x(x,z)\Big|_S\ -\    \mathcal{R}_z(x,z)\Big|_S \Big]
\end{split}
\end{align}
\textit{where}
\begin{align}\label{R(x,q)}
\mathcal{R}(x,z)=  \omega^*\Gamma_\eta(\mathbf{x}, \mathbf{q})\ G\xi+ \omega^*\Big(\Gamma_{q_1}(\mathbf{x},\mathbf{q})\dot{q}_1  + \Gamma_{q_2}(\mathbf{x},\mathbf{q})\dot{q}_2\Big)
\end{align}
\begin{align}\label{82}
 \Gamma(\mathbf{x},\mathbf{q})\ =&\   \frac{1}{4\pi}\ln \Big((x-q_1)^2\ +\ (z-q_2)^2\Big) \nonumber\\
&-\
   \sinh\big((z+h)D\big)\left[: \sinh\big((\eta+h)D\big):\right]^{-1} \Gamma_0(x,\eta,\mathbf{q}) \nonumber\\
&+\ :\sinh\big((z-\eta)D\big):\left[:\sinh\big((\eta+h)D\big):\right]^{-1} \Gamma_0(x,-h,\mathbf{q})\end{align}
\begin{align}\label{83}
\psi_x\Big|_S\ =&\   \omega^*  \Gamma_x(\mathbf{x},\mathbf{q})\Big|_S
-   :\sinh\big((\eta+h)D\big): D \left[:\cosh\big((\eta+h)D\big):\right]^{-1}\xi
\end{align}
\begin{align}\label{84}
\psi_z\Big|_S\ =&\ \omega^*  \Gamma_z(\mathbf{x},\mathbf{q})\Big|_S
+i: \cosh\big((\eta+h)D\big):D\left[:\cosh\big((\eta+h)D\big):\right]^{-1}\xi\end{align}

\section{Small-amplitude long-waves approximation}

For the propagation of  long waves, we will study the equations \eqref{eqs_point_v} under the assumption that their wavelength $\lambda$  is much bigger than the total depth of the undisturbed water $h$.
We denote by $\delta=\frac{h}{\lambda}$. Thus, for the wave number $k=\frac{2\pi}{\lambda}$ we have $k=\mathcal{O}(\delta)$. Because the operator $D$ has $k$ as eigenvalue (i.e.  $D{\rm e}^{ikx}=k{\rm e}^{ikx}$), we have $D=\mathcal{O}(\delta)$. By Taylor series expansion we have
\begin{align}
 \sinh  \big((z+h)D\big) \ =\   (z+h)D\ +\
\frac{(z+h)^3 D^3}{6}\ +\ \mathcal{O}(\delta^5),
\end{align}
\begin{align}
 \cosh \big((z+h)D\big) \ =\   1\ +\
\frac{(z+h)^2 D^2}{2}\ +\ \mathcal{O}(\delta^4).
\end{align}
In our analysis, we also consider that the waves have a small amplitude of order $\epsilon$.
For $\eta\in\mathcal{O}(\epsilon)$,  we make the  approximations
\begin{align}
\sinh\big((\eta+h)D\big)\ =&\ \sinh(\eta\, D)\cosh(h\, D)+\cosh(\eta\, D)\sinh(h\, D)\nonumber\\
  \approx & \ \eta\, D\cosh(h\, D)\ +\ \sinh(h\, D)
\end{align}
and
\begin{align}
\cosh\big((\eta+h)D\big)\ =&\ \cosh(\eta\, D)\cosh(h\, D)+\sinh(\eta\, D)\sinh(h\, D)\nonumber\\
  \approx & \ \cosh(h\, D)\ +\ \eta\, D\sinh(h\, D).
\end{align}
Therefore,
\begin{align}
\left[: \sinh\big((\eta+h)D\big):\right]^{-1}\approx & \left[:\big(\eta\, D \cosh(h\, D)+\sinh(h\, D):\right]^{-1}\nonumber\\
&= \left[:\big[\eta\, D \coth(h\, D)+1 \big]\sinh(h\, D):\right]^{-1}\nonumber\\
&=[  \sinh(h\, D)]^{-1}\left[:\big(\eta\, D \coth(h\, D)+1 \big):\right]^{-1}\nonumber\\
&=[  \sinh(h\, D)]^{-1}\left[:\big(1+\eta\, D \coth(h\, D) \big):\right]^{-1}\nonumber\\
&\approx\frac{1}{\sinh(h\, D)}\big(1-\eta\, D \coth(h\, D)\big)
\end{align}
and
\begin{align}
\left[: \cosh\big((\eta+h)D\big):\right]^{-1}\approx & \left[:\big(\eta\, D \sinh(h\, D)+\cosh(h\, D):\right]^{-1}\nonumber\\
&= \left[:\big[\eta\, D \tanh(h\, D)+1 \big]\cosh(h\, D):\right]^{-1}\nonumber\\
&=[  \cosh(h\, D)]^{-1}\left[:\big(\eta\, D \tanh(h\, D)+1 \big):\right]^{-1}\nonumber\\
&\approx\frac{1}{\cosh(h\, D)}\big(1-\eta\, D \tanh(h\, D)\big).
\end{align}
Further on,  by ignoring the terms of order $\mathcal{O}(\eta^2)$ , we get
\begin{align}
: \sinh & \big((\eta+h)D\big): D\, \left[: \cosh\big((\eta+h)D\big):\right]^{-1}\nonumber\\
\approx&\
\big(\sinh(h\, D)+\eta\, D\cosh(h\, D) \big)
D\, \frac{1}{\cosh(h\, D)}\big(1-\eta\, D \tanh(h\, D)\big)\nonumber\\
=& D \tanh(hD)- D\tanh(hD)\,  \eta\,  D\tanh(hD)+ \eta D^2\end{align}
and thus,  the expression \eqref{G2} of the Dirichlet-Neumann operator becomes
\begin{align}\label{G3}
G \xi\ =&\  -\eta_x\xi_x+  \big(D \tanh(hD)-  D\tanh(hD)\,  \eta\,  D \tanh(hD) +  \eta D^2\big) \xi+\mathcal{O}(\eta^2)\nonumber\\
\stackrel{\eqref{D}}{=}&\   (D\eta)\cdot( D\xi)+  D \tanh(hD) \xi +\eta D^2\xi- D\tanh(hD)\,  \eta\,  D \tanh(hD)  \xi+\mathcal{O}(\eta^2)\nonumber\\
=& D \tanh(hD) \xi + \big(D \eta D - D\tanh(hD)\,  \eta\,  D \tanh(hD)\big)\xi+\mathcal{O}(\eta^2)
\end{align}
We recognize here the first two  terms $$G^{(0)}\xi\ =\ D \tanh(hD) \xi$$ and $$G^{(1)}\xi\ =\ \big(D \eta D - D\tanh(hD)\,  \eta\,  D \tanh(hD)\big)\xi $$ from the Taylor series expansion of the Dirichlet-Neumann operator  \cite{Craig, Craig&Groves, CS}.\\
Next, by considering the expansion
\begin{align}
\tanh(hD)\ =\  hD\ -\  \frac{1}{3}h^3D^3\ +\  \mathcal{O}((hD)^5),
\end{align}
we get
\begin{align}
G \xi  \ =\  D h D\xi\ -\ \frac{1}{3}Dh^3D^3 \xi\ +\ D\eta D \xi\ +\ \mathcal{O}(\eta^2, \eta (hD)^2,  (hD)^4),
\end{align}
that is,
\begin{align}\label{Gapprox}
G\xi\ =\  -  h\xi_{xx}\ -\ \frac{1}{3}h^3\xi_{xxxx}\ -\ (\eta\xi_x)_x\ +\ \mathcal{O}(\epsilon^2, \epsilon\delta^2,\delta^4).
\end{align}
 We approximate the terms from the first and the second equations in \eqref{eqs_point_v}  by the following expressions
\begin{align}
 \cosh & \big((z+h)D\big) D\left[: \cosh\big((\eta+h)D\big):\right]^{-1}\nonumber\\
\approx&\  \Big(1+
\frac{(z+h)^2 D^2}{2}\Big)D \frac{1}{\cosh(h\, D)}\big(1-\eta\, D \tanh(h\, D)\big)\nonumber\\
\approx &\  \Big(1+
\frac{(z+h)^2 D^2}{2}\Big) D\Big(1-\frac{h^2D^2}{2}\Big)\Big(1-\eta\, D\Big(hD-\frac{h^3D^3}{6}\Big) \Big)\nonumber\\
\approx &\  D\ +\ \mathcal{O}(D^3),
\end{align}
\begin{align}
 \sinh & \big((z+h)D\big) D\left[: \cosh\big((\eta+h)D\big):\right]^{-1}\nonumber\\
\approx&\   (z+h)D\Big(1+
\frac{(z+h)^2 D^2}{6}\Big)D \frac{1}{\cosh(h\, D)}\big(1-\eta\, D \tanh(h\, D)\big)\nonumber\\
\approx &\   (z+h)D\Big(1+
\frac{(z+h)^2 D^2}{6}\Big) D\Big(1-\frac{h^2D^2}{2}\Big)\Big(1-\eta\, D\Big(hD-\frac{h^3D^3}{6}\Big) \Big)\nonumber\\
\approx &\ (z+h) D^2\ +\ \mathcal{O}(D^3).
\end{align}
Therefore,  these two equations become
\begin{align}\label{99}
 \begin{split}
\dot{q_1}\
=&\   \xi_x (q_1) \\
 \dot{q_2}\
 =&-  (q_2+h)\xi_{xx}(q_1),
\end{split}
\end{align}
that is, the point vortex moves only   because of the presence of the free surface and of the flat bottom.
The leading-order equations for the system  \eqref{99} are
\begin{align}
 \begin{split}
\dot{q_1}\
=&\ 0\\
 \dot{q_2}\
 =&\ 0
\end{split}
\hspace{1cm}
\implies \begin{split}
 \textrm{the point vortex stays at rest} \\
q_1(t)=q_1(0),\,  q_2(t)=q_2(0).
\end{split}
\end{align}

\noindent Similarly, we find for \eqref{83} and \eqref{84}
\begin{align}\label{100}
\begin{split}
\psi_x\Big|_S\ =&\   \omega^*  \Gamma_x(\mathbf{x},\mathbf{q})\Big|_S
+  (\eta+h)\xi_{xx}(x)\\
\psi_z\Big|_S\ =&\ \omega^*  \Gamma_z(\mathbf{x},\mathbf{q})\Big|_S
+\xi_x(x).
\end{split}
\end{align}
In  small-amplitude, long-wavelength regimes, the appropriate Green function for our point-vortex problem is the Green function on the infinite strip
$\{-\infty<x<\infty,\,   -h\ <\ z\ <0\}$. Alternative expressions of the  Green function on an infinite strip, constructed by different methods - the method of images, the method of conformal mapping or the method of eigenfunction expression - are available in standard books on partial differential equations, see, for example,  \cite{M&Pulvirenti, Melnikov} and the references therein. By the method of images one  obtains an infinite product representation of the  Green function
for the Dirichlet problem on an infinite strip  \cite{Melnikov}
\begin{align}
\mathit{\Gamma}(\mathbf{x},\mathbf{x}_0)=
\frac{1}{4\pi}\ln \prod_{n=-\infty}^{n=\infty}\frac{(x-x_0)^2+(z-z_0-2\,n\,h)^2}{(x-x_0)^2+(z+z_0+2\, n\, h)^2}.\end{align}
By the method of conformal mapping,   the following  Green function  is obtained \cite{Melnikov}
\begin{align}\label{Gamma-M}
\mathit{\Gamma}(\mathbf{x},\mathbf{x}_0)=\frac{1}{4\pi}\ln\frac{\cosh\frac{\pi}{h}(x-x_0)-\cos\frac{\pi}{h}(z-z_0)}{
\cosh\frac{\pi}{h}(x-x_0)-\cos\frac{\pi}{h}(z+z_0)}.
\end{align}
By taking $\eta=0$ in \eqref{61}-\eqref{62}, we obtain an alternative form of  the Green function on the infinite strip $\{-\infty<x<\infty,\,   -h\ <\ z\ <0\}$, that is,
\begin{align}\label{102}
 \mathit{\Gamma}(\mathbf{x},\mathbf{x}_0)\ =& \frac{1}{4\pi}\Big\{ \ln \Big((x-x_0)^2\ +\ (z-z_0)^2\Big) \nonumber\\
&-\
   \sinh\big((z+h)D\big) \left[\sinh\big(hD\big)\right]^{-1} \ln \Big((x-x_0)^2\ +\   z_0^2\Big) \nonumber\\
&+\  \sinh\big(z\, D\big)\left[\sinh\big(hD\big)\right]^{-1} \ln \Big((x-x_0)^2\ +\ (h+z_0)^2\Big)\Big\}
\end{align}
For the following calculations,  we will choose the expression \eqref{Gamma-M} of the Green function
\begin{align}\label{Gamma-M-q}
\mathit{\Gamma}(\mathbf{x},\mathbf{q})=\frac{1}{4\pi}\ln\frac{\cosh\frac{\pi}{h}(x-q_1)-\cos\frac{\pi}{h}(z-q_2)}{
\cosh\frac{\pi}{h}(x-q_1)-\cos\frac{\pi}{h}(z+q_2)}.
\end{align}
We see that far away from the point vortex, that is, $| \mathbf{x}-\mathbf{q} | \rightarrow \infty$, we have $\mathit{\Gamma}(\mathbf{x},\mathbf{q}) \rightarrow 0$,  and close to the vortex, that is,  $\mathbf{x}\rightarrow\mathbf{q}$, the Green function diverges as a logarithm.
 By accepting, as in \cite{M&Pulvirenti, Meleshko}, that there is no self-interaction, the above formula and the derivatives
\begin{align}
\begin{split}
 \mathit{\Gamma}_x(\mathbf{x},\mathbf{q})=\frac{1}{2h}\frac{\sin\left(\frac{\pi}{h}z\right) \sin\left(\frac{\pi}{h}q_2\right)\sinh\frac{\pi}{h}(x-q_1)}{
\big[\cosh\frac{\pi}{h}(x-q_1)-\cos\frac{\pi}{h}(z-q_2)\big]\big[\cosh\frac{\pi}{h}(x-q_1)-\cos\frac{\pi}{h}(z+q_2)\big]}\\
 \mathit{\Gamma}_z(\mathbf{x},\mathbf{q})=\frac{1}{2h}\frac{
\sin\left(\frac{\pi}{h}q_2\right)
\big[\cos\left(\frac{\pi}{h}q_2\right)-\cosh\frac{\pi}{h}(x-q_1) \cos\left(\frac{\pi}{h}z\right)\big]
}{
\big[\cosh\frac{\pi}{h}(x-q_1)-\cos\frac{\pi}{h}(z-q_2)\big]\big[\cosh\frac{\pi}{h}(x-q_1)-\cos\frac{\pi}{h}(z+q_2)\big]}\label{Gamma_x&z}
\end{split}
\end{align}
 apply everywhere   $\mathbf{x}\neq \mathbf{q}$. \\
In the system \eqref{100},
we need
 the expressions of $\mathit{\Gamma}_x(\mathbf{x},\mathbf{q})\Big|_S$ and $\mathit{\Gamma}_z(\mathbf{x},\mathbf{q})\Big|_S$,
and further, in the fifth equation of  \eqref{eqs_point_v},  we need the expression
\begin{align}\label{111}
 \begin{split}
|\nabla\psi|^2_S=(\psi_x\Big|_S+\psi_z\Big|_S)^2\ =&\   (\omega^*)^2  \left(\Gamma_x^2(\mathbf{x},\mathbf{q})\Big|_S+\Gamma_z^2(\mathbf{x},\mathbf{q})\Big|_S\right)\\
&+  (\eta+h)^2\xi_{xx}^2(x)+\xi_x^2(x)\\
&+ 2\omega^*  (\eta+h)\xi_{xx}\Gamma_x(\mathbf{x},\mathbf{q})\Big|_S\\
&+2\omega^*\xi_x(x)\Gamma_z(\mathbf{x},\mathbf{q})\Big|_S.
\end{split}
\end{align}
For a  free surface $\eta$ of order $\epsilon$,    the functions  $\mathit{\Gamma}_x(\mathbf{x},\mathbf{q})\Big|_S=\Gamma_x(x,\eta,q_1,q_2)$ and $\mathit{\Gamma}_z(\mathbf{x},\mathbf{q})\Big|_S=\Gamma_z(x,\eta,q_1,q_2)$,  have the following  Taylor series expansions
$$
  \Gamma_x(x, \eta; q_1,q_2)=\Gamma_x(x, 0; q_1,q_2) +  \eta \Gamma_{xz}(x, 0; q_1,q_2) + ...=\Gamma_x(x, 0; q_1,q_2)+\mathcal{O}(\epsilon)
$$
$$
  \Gamma_z(x, \eta; q_1,q_2)=\Gamma_z(x, 0; q_1,q_2) + \eta \Gamma_{zz}(x, 0; q_1,q_2) + ...=\Gamma_z(x, 0; q_1,q_2) +\mathcal{O}(\epsilon)
$$
From \eqref{Gamma_x&z}, we get
\begin{align}\label{105}
\begin{split}
 \mathit{\Gamma}_x(x,0,q_1,q_2) &= 0\\
 \mathit{\Gamma}_z(x,0,q_1,q_2) &=- \frac{\sin\frac{\pi q_2}{h}}{
2h \left(\cosh \frac{\pi(x-q_1)}{h}-\cos \frac{\pi q_2}{h} \right)},
\end{split}
\end{align}
therefore,
\begin{align}\label{113}
\begin{split}
 \mathit{\Gamma}_x(\mathbf{x},\mathbf{q})\Big|_S & = \mathcal{O}(\epsilon)\\
 \mathit{\Gamma}_z(\mathbf{x},\mathbf{q})\Big|_S &= \mathit{\Gamma}_z(x,0,q_1,q_2)+\mathcal{O}(\epsilon).
\end{split}
\end{align}

The $x$-derivative of the velocity potential on the free surface does not get an extra factor of $\delta$, since the $\epsilon$ order of the velocity remains unchanged.  Thus,  $\xi_x\in \mathcal{O}(\epsilon)$  and $\xi_{xx}\in \mathcal{O}(\epsilon\delta)$.
We make the notation
\begin{align}\label{109}
\mathfrak{u}:=\xi_x \in \mathcal{O}(\epsilon)
\end{align}
thus,
\begin{align}\label{110}
\mathfrak{u}_x=\xi_{xx} \in \mathcal{O}(\epsilon\delta).
\end{align}
In what follows we consider the wave propagation regime with $\epsilon \sim \delta^2$ which usually leads to the Boussinesq and KdV propagation regimes.\\
The magnitude of the vortex strength in this regime depends on the vortex depth and some estimates are made in the Appendix.

\noindent In view of \eqref{113},  \eqref{110}, \eqref{122} or \eqref{123}, we make the approximation
\begin{align}
|\nabla\psi|^2_S\ \approx &\    (\omega^*)^2  \Gamma_z^2(x,0,q_1,q_2)+2\omega^*\xi_x(x)\Gamma_z(x,0,q_1,q_2)+\xi_x^2(x)\nonumber\\
&=\big[\omega^*\Gamma_z(x,0,q_1,q_2)+\xi_x(x)\big]^2.
\end{align}
 For the Green function \eqref{Gamma-M-q} the derivative with respect to $\eta$ is zero, thus,  the function
$\mathcal{R}(x,z)$ in \eqref{R(x,q)} have in this case the expression
\begin{align}\label{119}
\mathcal{R}(x,z)=&  \omega^*\Big(\Gamma_{q_1}(\mathbf{x},\mathbf{q})\dot{q}_1  + \Gamma_{q_2}(\mathbf{x},\mathbf{q})\dot{q}_2\Big)\nonumber\\
\stackrel{\eqref{99}, \eqref{109}}{=}&  \omega^*\Big(\Gamma_{q_1}(\mathbf{x},\mathbf{q})\mathfrak{u}(q_1) -(q_2+h) \Gamma_{q_2}(\mathbf{x},\mathbf{q})\mathfrak{u}_{x}(q_1)\Big).
\end{align}
The fifth equation in \eqref{eqs_point_v} becomes
\begin{align}
\xi_t  +&\frac{1}{2}\left[
\omega^*\Gamma_z(x,0,q_1,q_2)+\xi_x(x)\right]^2  + g\eta\nonumber\\
& -  \partial_x^{-1}\Big[\eta_x \mathcal{R}_x(x,z)\Big|_S\ -\    \mathcal{R}_z(x,z)\Big|_S \Big]=0,
\end{align}
equation that, using also \eqref{109},  we write  as
\begin{align}\label{116}
\mathfrak{u}_t  +&\left\{\frac{1}{2}\left[
\omega^*\Gamma_z(x,0,q_1,q_2)+\mathfrak{u}\right]^2\right\}_x  + g\eta_x\nonumber\\
& -  \eta_x \mathcal{R}_x(x,z)\Big|_S\ +\    \mathcal{R}_z(x,z)\Big|_S =0.
\end{align}
From \eqref{119},
\begin{align}
\begin{split}
\mathcal{R}_x(x,z)&=  \omega^*\Big(\Gamma_{xq_1}(\mathbf{x},\mathbf{q})\mathfrak{u}(q_1) -(q_2+h) \Gamma_{xq_2}(\mathbf{x},\mathbf{q})\mathfrak{u}_{x}(q_1)\Big)\\
\mathcal{R}_z(x,z)&=  \omega^*\Big(\Gamma_{zq_1}(\mathbf{x},\mathbf{q})\mathfrak{u}(q_1) -(q_2+h) \Gamma_{zq_2}(\mathbf{x},\mathbf{q})\mathfrak{u}_{x}(q_1)\Big).
\end{split}
\end{align}
For the Green function \eqref{Gamma-M-q}, the expressions of the  second derivatives with respect to $x$ and $z$, are obtained further from  \eqref{Gamma_x&z}. Then, we need to evaluate these expressions on the free surface. By writing the Taylor series expansions of these functions
for a free surface $\eta$ of order $\epsilon$,   we get after calculations
\begin{align}\label{123'}
\begin{split}
 \mathit{\Gamma}_{xq_1}(\mathbf{x},\mathbf{q})\Big|_S & = \mathcal{O}(\epsilon)\\
 \mathit{\Gamma}_{zq_1}(\mathbf{x},\mathbf{q})\Big|_S &= -\frac{\pi}{2h^2}\frac{\sin\frac{\pi q_2}{h}\sinh\frac{\pi(x-q_1)}{h}}{\left(\cosh\frac{\pi(x-q_1)}{h}-\cos\frac{\pi q_2}{h}\right)^2}+\mathcal{O}(\epsilon^2)\\
 \Gamma_{xq_2}(\mathbf{x},\mathbf{q})\Big|_S&=\mathcal{O}(\epsilon)\\
 \mathit{\Gamma}_{zq_2}(\mathbf{x},\mathbf{q})\Big|_S &=\frac{\pi}{2h^2}\frac{1-\cos\frac{\pi q_2}{h}\cosh\frac{\pi(x-q_1)}{h}}{\left(\cosh\frac{\pi(x-q_1)}{h}-\cos\frac{\pi q_2}{h}\right)^2} +\mathcal{O}(\epsilon^2).
\end{split}
\end{align}
\noindent In view of   \eqref{109}, \eqref{110}, \eqref{123'}, \eqref{122} or \eqref{123}, we make the approximation
\begin{align}\label{124}
\begin{split}
 -  &\eta_x \mathcal{R}_x(x,z)\Big|_S +    \mathcal{R}_z(x,z)\Big|_S \approx   - \frac{\omega^*\pi}{2h^2}\frac{\sin\frac{\pi q_2}{h}\sinh\frac{\pi(x-q_1)}{h}}{\left(\cosh\frac{\pi(x-q_1)}{h}-\cos\frac{\pi q_2}{h}\right)^2} \, \mathfrak{u}(q_1)
\\
&\hspace{4cm} - (q_2+h)\frac{\omega^*\pi}{2h^2}\frac{1-\cos\frac{\pi q_2}{h}\cosh\frac{\pi(x-q_1)}{h}}{\left(\cosh\frac{\pi(x-q_1)}{h}-\cos\frac{\pi q_2}{h}\right)^2}\mathfrak{u}_{x}(q_1) \\
%&=- \frac{\omega^*\pi\left[\sin\frac{\pi q_2}{h}\sinh\frac{\pi(x-q_1)}{h}\mathfrak{u}(q_1) +  (q_2\! +\! h)(1-\cos\frac{\pi q_2}{h}\cosh\frac{\pi(x-q_1)}{h})\mathfrak{u}_{x}(q_1)\right]}{2h^2\left(\cosh\frac{\pi(x-q_1)}{h}-\cos\frac{\pi q_2}{h}\right)^2}
\end{split}
\end{align}

\noindent Summing up, from  \eqref{Gapprox},  \eqref{99}, \eqref{105}, \eqref{116} and \eqref{109}, \eqref{124},  \textit{we obtain in the regime $\epsilon \sim \delta^2$, the following system for the dynamics of the point vortex and the evolution of the free surface variables}
\begin{align}\label{system-final}
 \begin{split}
& \dot{q_1}\
=  \mathfrak{u} (q_1(t),t)\\
& \dot{q_2}\
 = - (q_2+h)\mathfrak{u}_{x}(q_1)\\
& \eta_t\ + h\mathfrak{u}_{x}\ +\   \frac{1}{3}h^3\mathfrak{u}_{xxx}\ +\  (\eta\mathfrak{u})_x=0\\
& \mathfrak{u}_t +g\eta_x+ \left\{\frac{1}{2}\left[\mathfrak{u}-\omega^*\frac{\sin\frac{\pi q_2}{h}}{
2h \left(\cosh \frac{\pi(x-q_1)}{h}-\cos \frac{\pi q_2}{h} \right)}\right]^2 \right \}_x\\
&\hspace{2cm}- \frac{\omega^*\pi}{2h^2}\left\{ \frac{\sin\frac{\pi q_2}{h}\sinh\frac{\pi(x-q_1)}{h}}{\left(\cosh\frac{\pi(x-q_1)}{h}-\cos\frac{\pi q_2}{h}\right)^2} \, \mathfrak{u}(q_1)\right.\\
&\hspace{3cm}\left.+(q_2+h)\frac{1-\cos\frac{\pi q_2}{h}\cosh\frac{\pi(x-q_1)}{h}}{\left(\cosh\frac{\pi(x-q_1)}{h}-\cos\frac{\pi q_2}{h}\right)^2}\mathfrak{u}_{x}(q_1) \right\}=0.
\end{split}
\end{align}

\section{Derivation of the perturbed KdV equation}

For $\eta,\, \mathfrak{u},\, \omega^*\in \mathcal{O}(\epsilon)$ (see the appendix),  $\eta_t, \, \mathfrak{u}_x\in \mathcal{O}(\epsilon\delta)$, $\mathfrak{u}_{xxx},\, (\eta\mathfrak{u})_x\in \mathcal{O}(\epsilon^2\delta)$,  let us introduce for the last two equations in \eqref{system-final} the scaling
\begin{align}\label{CV1}
 \begin{split}
& \eta_t + h\mathfrak{u}_{x} +  \epsilon \frac{1}{3}h^3\mathfrak{u}_{xxx}\ + \epsilon (\eta\mathfrak{u})_x=0\\
& \mathfrak{u}_t +g\eta_x\
+ \epsilon  \left\{\frac{1}{2}\left[\mathfrak{u}-\omega^*\frac{\sin\frac{\pi q_2}{h}}{
2h \left(\cosh \frac{\pi(x-q_1)}{h}-\cos \frac{\pi q_2}{h} \right)}\right]^2 \right \}_x\\
&\hspace{1.7cm}-\epsilon\, \frac{\omega^*\pi}{2h^2}\left\{\frac{\sin\frac{\pi q_2}{h}\sinh\frac{\pi(x-q_1)}{h}}{\left(\cosh\frac{\pi(x-q_1)}{h}-\cos\frac{\pi q_2}{h}\right)^2} \, \mathfrak{u}(q_1)\right.\\
&\hspace{3cm}\left. + (q_2+h)\frac{1-\cos\frac{\pi q_2}{h}\cosh\frac{\pi(x-q_1)}{h}}{\left(\cosh\frac{\pi(x-q_1)}{h}-\cos\frac{\pi q_2}{h}\right)^2}\mathfrak{u}_{x}(q_1) \right\}=0.
\end{split}
\end{align}
Next, for describing running waves in one $x$-direction  with the speed
\begin{align}\label{speed}
c^2= gh,
\end{align}
 we move  from the  $(x,t)$-frame to the $(X,T)$-frame, where the characteristic variable $X$ and the slow time $T$ are defined by
\begin{align}\label{XT}
X = x-ct, \quad T=\epsilon t.
\end{align}
 In terms of these variables, we have
\begin{align}\label{CoV}
 \partial_x= \partial_{X} ,\quad \partial_t=\epsilon \partial_T - c \partial_{X},
\end{align}
and  the equations \eqref{CV1} become
\begin{align}\label{CV2}
 \begin{split}
&\epsilon \eta_T  -c\eta_X + h\mathfrak{u}_{X}   +  \epsilon \frac{1}{3}h^3\mathfrak{u}_{X X X}\ + \epsilon (\eta\mathfrak{u})_{X}=0\\
& \epsilon  \mathfrak{u}_T - c \mathfrak{u}_{X} +g\eta_{X}+ \epsilon
 \left\{\frac{1}{2}\left[\mathfrak{u}-\omega^*\frac{\sin\frac{\pi q_2}{h}}{
2h \left(\cosh \frac{\pi(\tilde{x}-q_1)}{h}-\cos \frac{\pi q_2}{h} \right)}\right]^2 \right \}_X\\
&\hspace{1.7cm}
-\epsilon\frac{\omega^*\pi}{2h^2}\left\{ \frac{\sin\frac{\pi q_2}{h}\sinh\frac{\pi(\tilde{x}-q_1)}{h}}{\left(\cosh\frac{\pi(\tilde{x}-q_1)}{h}-\cos\frac{\pi q_2}{h}\right)^2} \, \mathfrak{u}(q_1)\right.\\
&\hspace{3cm}\left.+(q_2+h)\frac{1-\cos\frac{\pi q_2}{h}\cosh\frac{\pi(\tilde{x}-q_1)}{h}}{\left(\cosh\frac{\pi(\tilde{x}-q_1)}{h}-\cos\frac{\pi q_2}{h}\right)^2}\mathfrak{u}_{X}(q_1)\right\} =0,
\end{split}
\end{align}
with $\tilde{x}=x(X, T).$ From the first equation, we get

\begin{equation}\label{CV3}
   \eta_{X}=  \frac{h}{c}\mathfrak{u}_{X} + \epsilon \frac{1}{c} \eta_T+ \epsilon \frac{h^3}{3c}\mathfrak{u}_{X X X}\ + \epsilon \frac{1}{c}(\eta\mathfrak{u})_{X},
\end{equation}
which yields \begin{align}\label{132}
\eta = \frac{h}{c}\mathfrak{u} + \mathcal{O}(\epsilon).
\end{align}
 We substitute \eqref{CV3} into the second equation of \eqref{CV2}. Due to \eqref{speed},  the leading order terms with $\mathfrak{u}_{X} $ cancel  out  and we get
\begin{equation}\label{CV5}
\begin{split}
  &    \mathfrak{u}_T  +  \frac{g}{c} \eta_T+  \frac{g h^3}{3c}\mathfrak{u}_{X X X}\ + \frac{g}{c}(\eta\mathfrak{u})_{X}\\
& \hspace{2.37cm}+
 \left\{\frac{1}{2}\left[\mathfrak{u}-\omega^*\frac{\sin\frac{\pi q_2}{h}}{
2h \left(\cosh \frac{\pi(\tilde{x}-q_1)}{h}-\cos \frac{\pi q_2}{h} \right)}\right]^2 \right \}_X\\
&\hspace{1.7cm}-\frac{\omega^*\pi}{2h^2}\left\{ \frac{\sin\frac{\pi q_2}{h}\sinh\frac{\pi(\tilde{x}-q_1)}{h}}{\left(\cosh\frac{\pi(\tilde{x}-q_1)}{h}-\cos\frac{\pi q_2}{h}\right)^2} \, \mathfrak{u}(q_1)\right.\\
&\hspace{3cm}\left.+ (q_2+h)\frac{1-\cos\frac{\pi q_2}{h}\cosh\frac{\pi(\tilde{x}-q_1)}{h}}{\left(\cosh\frac{\pi(\tilde{x}-q_1)}{h}-\cos\frac{\pi q_2}{h}\right)^2}\mathfrak{u}_{X}(q_1) \right\} =0,
\end{split}
\end{equation}
where all terms are of the same order.  With \eqref{132} and \eqref{speed} in view, by keeping only the leading order terms in \eqref{CV5}, we obtain
\begin{equation}\label{CV8}
\begin{split}
&     2\mathfrak{u}_T  +  \frac{c h^2}{3}\mathfrak{u}_{XXX}\ +
 \left\{\frac{1}{2}\left[\mathfrak{u}-\omega^*\frac{\sin\frac{\pi q_2}{h}}{
2h \left(\cosh \frac{\pi(\tilde{x}-q_1)}{h}-\cos \frac{\pi q_2}{h} \right)}\right]^2+\mathfrak{u}^2 \right \}_X\\
&\hspace{1.7cm}- \frac{\omega^*\pi}{2h^2}\left\{\frac{\sin\frac{\pi q_2}{h}\sinh\frac{\pi(\tilde{x}-q_1)}{h}}{\left(\cosh\frac{\pi(\tilde{x}-q_1)}{h}-\cos\frac{\pi q_2}{h}\right)^2} \, \mathfrak{u}(q_1)\right.\\
&\hspace{3cm}\left.+ (q_2+h)\frac{1-\cos\frac{\pi q_2}{h}\cosh\frac{\pi(\tilde{x}-q_1)}{h}}{\left(\cosh\frac{\pi(\tilde{x}-q_1)}{h}-\cos\frac{\pi q_2}{h}\right)^2}\mathfrak{u}_{X}(q_1) \right\} =0.
\end{split}
\end{equation}
This is a perturbed KdV equation.\\
 If $\omega^*=0,$ we have the standard KdV equation,
\begin{equation}\label{CV9}
     \mathfrak{u}_T  +  \frac{c h^2}{6}\mathfrak{u}_{X X X}
   +    \frac{3}{4}(\mathfrak{u}^2)_{X}=0.
\end{equation}
We return  now to the original variables \eqref{XT}-\eqref{CoV} and the equation \eqref{CV8} becomes
\begin{align}\label{CV10}
\begin{split}
     \mathfrak{u}_t& +c   \mathfrak{u}_x+\epsilon \frac{c h^2}{6}\mathfrak{u}_{xxx}\\
&  \hspace{0.6cm}+   \epsilon
\left\{\frac{1}{4}\left[\mathfrak{u}-\omega^*\frac{\sin \frac{\pi q_2}{h}}{2h \left(\cosh \frac{\pi(x-q_1)}{h}-\cos \frac{\pi q_2}{h} \right)}\right]^2 +\frac{1}{2} \mathfrak{u}^2 \right\}_{x}\\
&  \hspace{0.7cm}-\epsilon\frac{\omega^*\pi}{4h^2}\left\{\frac{\sin\frac{\pi q_2}{h}\sinh\frac{\pi(x-q_1)}{h}}{\left(\cosh\frac{\pi(x-q_1)}{h}-\cos\frac{\pi q_2}{h}\right)^2} \, \mathfrak{u}(q_1)\right.\\
&\hspace{2cm}\left.+ (q_2+h)\frac{1-\cos\frac{\pi q_2}{h}\cosh\frac{\pi(x-q_1)}{h}}{\left(\cosh\frac{\pi(x-q_1)}{h}-\cos\frac{\pi q_2}{h}\right)^2}\mathfrak{u}_{x}(q_1) \right\}
 =0.
\end{split}
\end{align}
Without the scaling written explicitly, \textit{we get the following   perturbed KdV equation for $\mathfrak{u}(x,t)$ coupled with the dynamic equations for the point vortex $\left(q_1(t),q_2(t)\right)$}
\begin{align}\label{KdV system q}
 \begin{split}
&    \mathfrak{u}_t  +c   \mathfrak{u}_x+ \frac{c h^2}{6}\mathfrak{u}_{xxx}\\
&\hspace{0.7cm}+
\left\{\frac{1}{4}\left[\mathfrak{u}-\omega^*\frac{\sin \frac{\pi q_2}{h}}{2h \left(\cosh \frac{\pi(x-q_1)}{h}-\cos \frac{\pi q_2}{h} \right)}\right]^2 +\frac{1}{2} \mathfrak{u}^2 \right\}_{x}\\
&\hspace{0.7cm}-\frac{\omega^*\pi}{4h^2}\left\{\frac{\sin\frac{\pi q_2}{h}\sinh\frac{\pi(x-q_1)}{h}}{\left(\cosh\frac{\pi(x-q_1)}{h}-\cos\frac{\pi q_2}{h}\right)^2} \, \mathfrak{u}(q_1)\right.\\
&\hspace{2cm}\left.+ (q_2+h)\frac{1-\cos\frac{\pi q_2}{h}\cosh\frac{\pi(x-q_1)}{h}}{\left(\cosh\frac{\pi(x-q_1)}{h}-\cos\frac{\pi q_2}{h}\right)^2}\mathfrak{u}_{x}(q_1) \right\}  =0\\
& \dot{q_1}\
=  \mathfrak{u} (q_1(t),t)\\
&
 \dot{q_2}\
 =-  (q_2+h)\mathfrak{u}_{x}(q_1(t),t) .
\end{split}
\end{align}

\section{Vortex motion influenced by an unperturbed KdV soliton }

The  coupled system \eqref{KdV system q} is quite  complicated and deserves a proper numerical investigation. However, some features could be revealed
by the simplification that the soliton on the surface travels unperturbed above the vortex. The exact KdV soliton has the expression (see, for example, \cite{Cullen&Ivanov})
\begin{equation}\label{KdVsol}
\mathfrak{u}(x,t)=\frac{4 c h^2 K^2}{3 \cosh^2\left[K\left(x-x_0-(1+\frac{2}{3}h^2K^2)\, ct\right)\right]},
\end{equation} where $x_0$ and $K$, called the soliton parameters, are constants. Then, the last two equations in \eqref{KdV system q} can be written as
\begin{align}\label{qKdVsoliton}
 \begin{split}
& \dot{q_1}= \frac{4 c h^2 K^2}{3 \cosh^2\left[K\left(q_1-x_0-(1+\frac{2}{3} h^2K^2)\, ct\right)\right]}\\
&
 \dot{q_2} =  (q_2+h)\frac{8 c h^2 K^3  \sinh \left[K\left(q_1-x_0-(1+\frac{2}{3} h^2K^2)\, ct\right)\right] }{3 \cosh^3\left[K\left(
q_1-x_0-(1+\frac{2}{3}h^2K^2)\, ct\right)\right]}.
\end{split}
\end{align}
Without loss of generality, we can take $x_0=0,$  which corresponds to a change of variables $q_1 \rightarrow q_1 -x_0$. Thus, the trajectories of the point vortex are described by
\begin{align}\label{qKdVsoliton1}
 \begin{split}
& \dot{q_1}= \frac{4 c h^2 K^2}{3 \cosh^2\left[K\left(q_1-(1+\frac{2}{3} h^2K^2)\, ct\right)\right]}\\
&
 \dot{q_2} =  (q_2+h)\frac{8 c h^2 K^3  \sinh \left[K\left(q_1-(1+\frac{2}{3} h^2K^2)\, ct\right)\right] }{3 \cosh^3\left[K\left(q_1-(1+\frac{2}{3} h^2K^2)\, ct\right)\right]}.
\end{split}
\end{align}
The first equation of the above system can be integrated and the conservation law implicitly provides the solution,
\begin{align}\label{cl}
q_1(t)+\frac{2h }{\sqrt{3-2h^2K^2}  } \textrm{arctan}\left[ \frac{2hK \tanh \left[K\left(q_1(t)-(1+\frac{2}{3}h^2K^2)\, ct\right)\right]    }{\sqrt{3-2h^2K^2}} \right]=\text{const.}
\end{align}
for $3-2h^2K^2>0$. According to our assumptions, $K^2h^2$ is of order $\epsilon,$ thus $3-2h^2K^2>0$ indeed. Moreover, since $\dot{q_1}$ is of order $\epsilon,$ and $(1+\frac{2}{3}h^2K^2)\, c$ is of order 1, it follows that $q_1(t)-(1+\frac{2}{3}h^2K^2)\, ct \to \mp \infty$ as $t \to \pm \infty.$ Therefore, if we go to the limits $t\to \pm \infty$ in \eqref{cl}, and take into account that at these limits $\tanh\rightarrow \mp 1$, we get
\begin{align}
q_1(\infty)-q_1(-\infty)= \frac{4h}{\sqrt{3-2h^2K^2}} \textrm{arctan}\left[\frac{2hK}{\sqrt{3-2h^2K^2}}\right],
\end{align}
that is, the trajectories are bounded, the range of change being the above finite value.

By using Maple, we plot in  Fig. \ref{fig:q1},   Fig. \ref{fig:q2} and  Fig. \ref{fig:q1vsq2}, the graph of $q_1$ versus $t$, the graph of $q_2$ versus $t$, and the graph of $q_1$ versus $q_2$, that is, the point vortex trajectories, respectively.  In all the pictures, we take $h=10$m, $g=9.81$m/s$^2$ and  the initial condition  $q_1(0)=0$. The initial condition for $q_2$   takes different values:  $q_2(0)=-0.1$ (green), $q_2(0)=-6$ (red), $q_2(0)=-9.9$ (blue). For  the pictures in  $(a)$  we consider  $K=0.1$m$^{-1}$, that is, a stronger soliton compare to the  pictures in $(b)$ drawn for $K=0.05$m$^{-1}$.
\begin{figure}[h!]
\begin{center}
(a)\fbox{\includegraphics[totalheight=0.4\textheight]{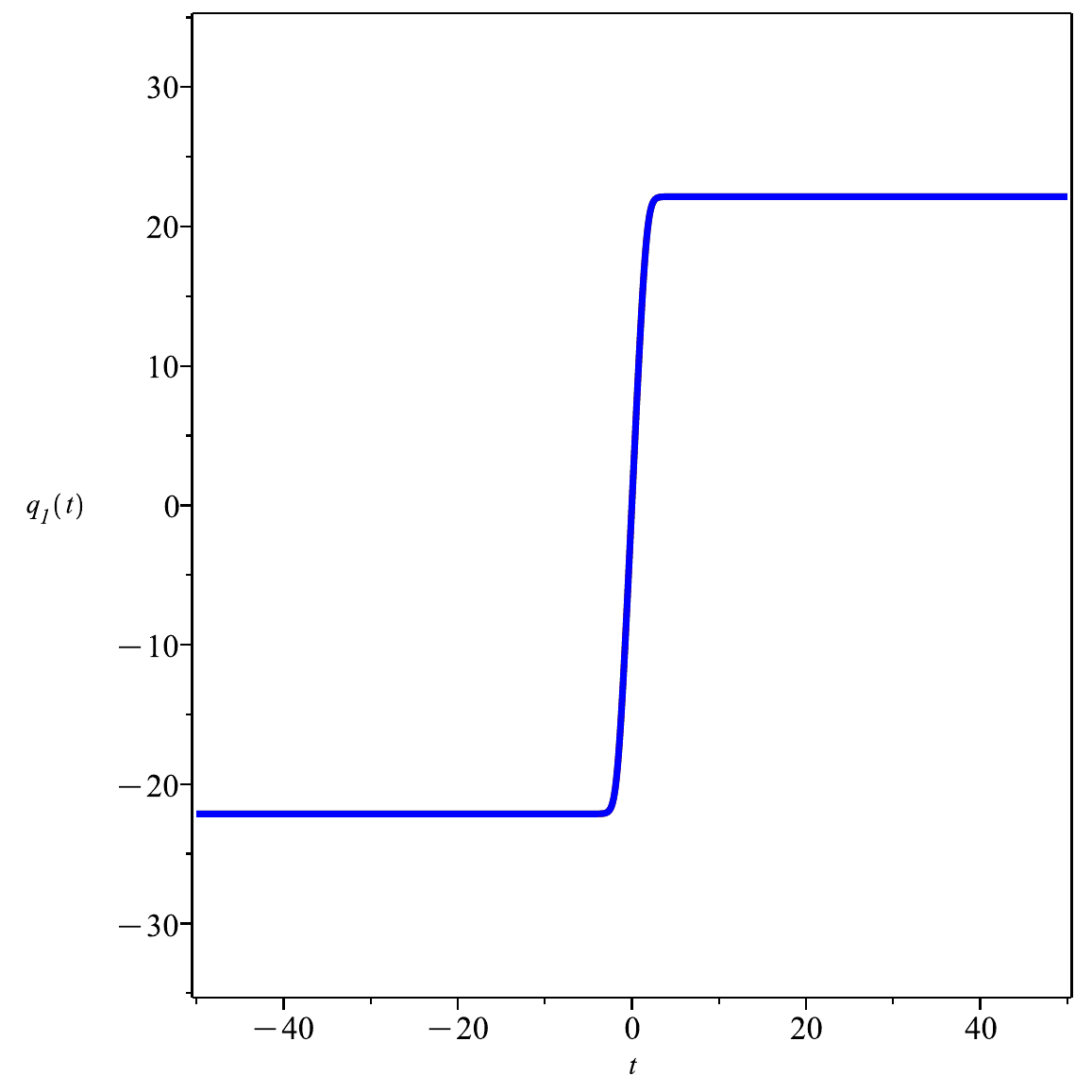}}
(b)\fbox{\includegraphics[totalheight=0.4\textheight]{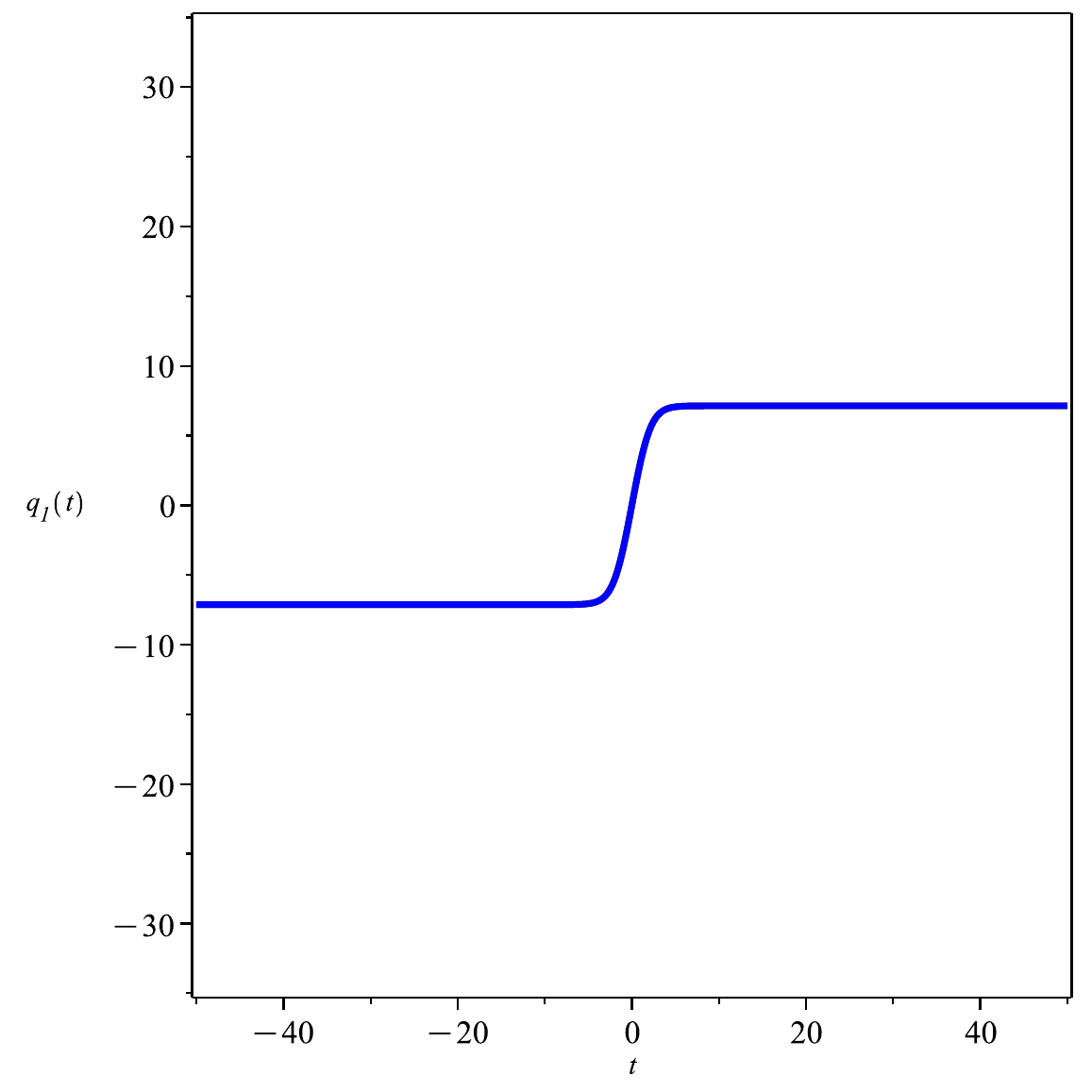}}
\caption{The evolution in time of $q_1$, with the initial condition $q_1(0)=0$ and  $h=10$m, in the cases: (a) $K=0.1$m$^{-1}$; (b) $K=0.05$m$^{-1}.$}
\label{fig:q1}
\end{center}
\end{figure}
\begin{figure}[h!]
\begin{center}
(a)\fbox{\includegraphics[totalheight=0.4\textheight]{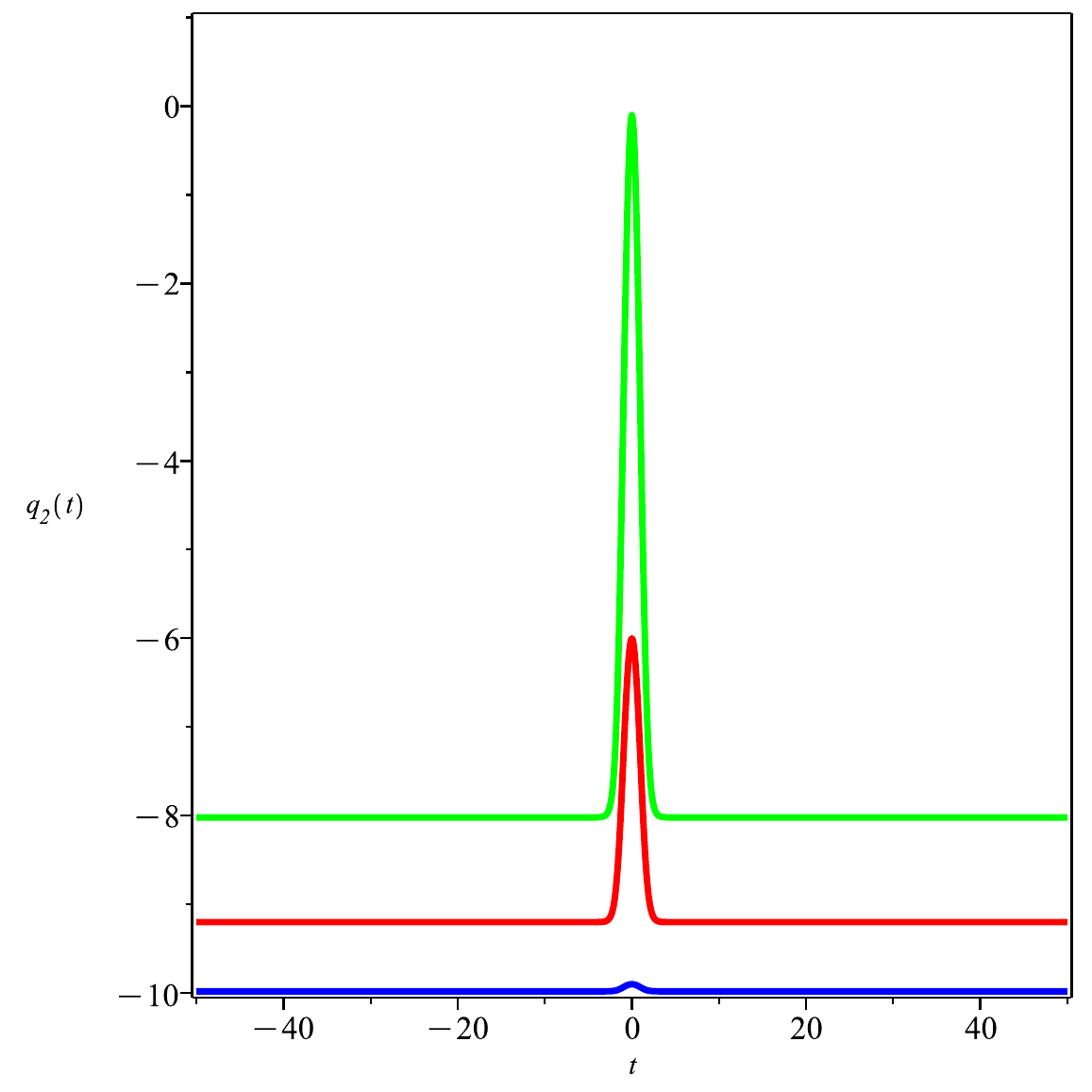}}
(b)\fbox{\includegraphics[totalheight=0.4\textheight]{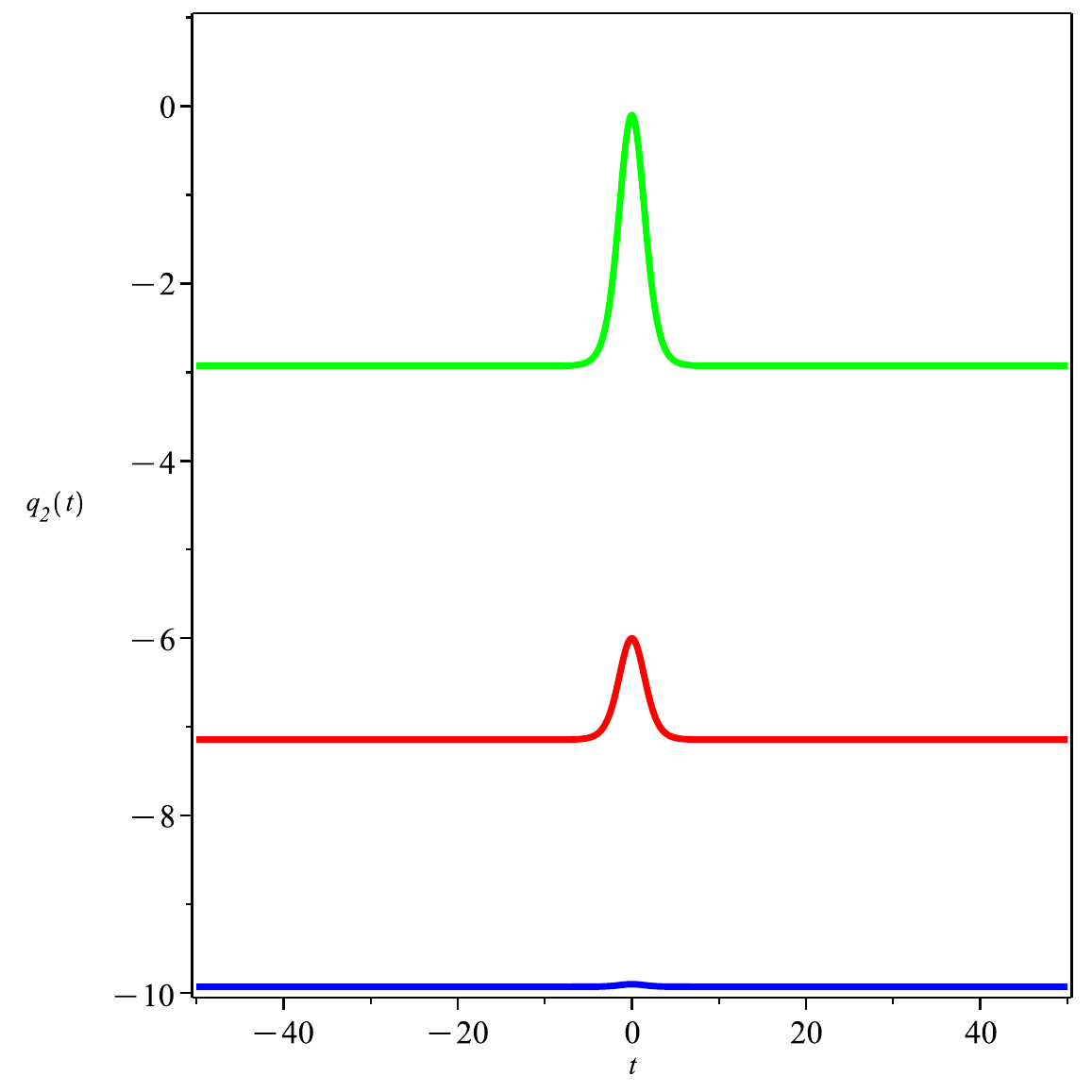}}
\caption{
The evolution in time of $q_2$, with the initial condition $q_2(0) = -0.1$ (up, green),
$q_2(0) = -6$ (in the middle, red), $q_2(0) = -9.9$ (down, blue) and $h=10$m, in the cases: (a) $K=0.1$m$^{-1}$; (b) $K=0.05$m$^{-1}.$}
\label{fig:q2}
\end{center}
\end{figure}
\begin{figure}[h!]
\begin{center}
(a)\fbox{\includegraphics[totalheight=0.4\textheight]{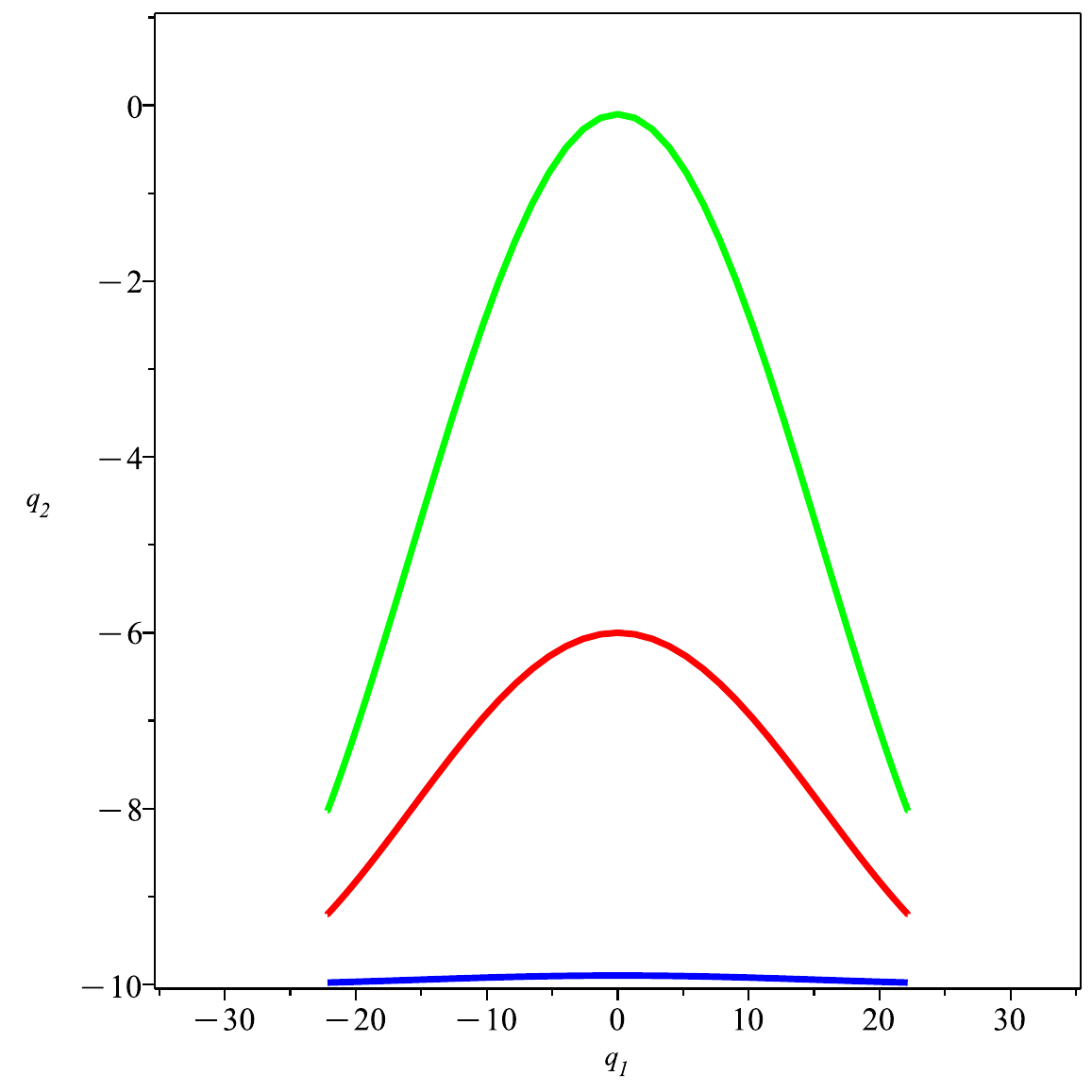}}
(b)\fbox{\includegraphics[totalheight=0.4\textheight]{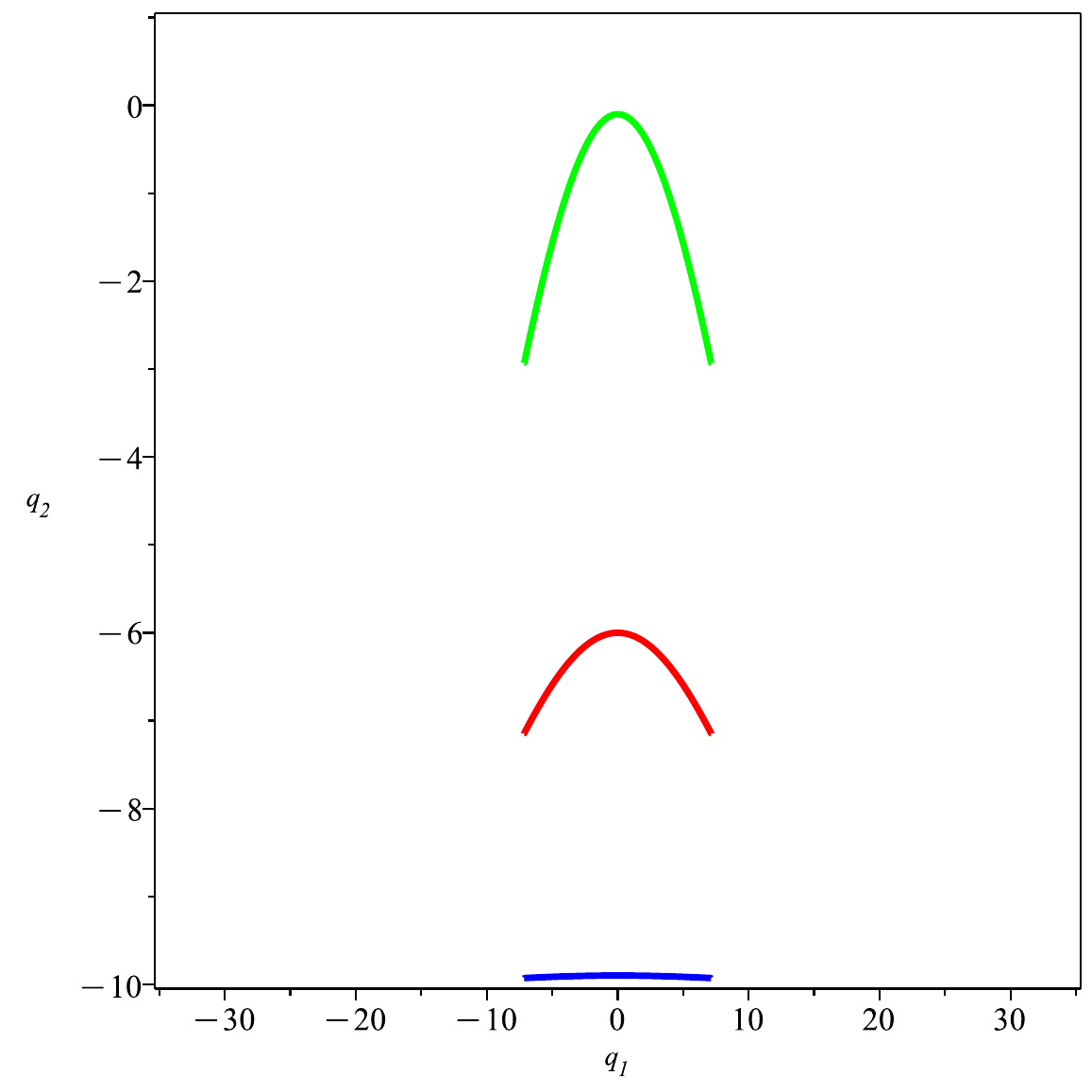}}
\caption{Vortex trajectories at different depths, with the  initial conditions  $q_1(0) = 0$ in all cases, $q_2(0) = -0.1$ (up, green),
$q_2(0) = -6$ (in the middle, red), $q_2(0) = -9.9$ (down, blue) and $h=10$m,
 in the cases: (a) $K=0.1$m$^{-1}$; (b) $K=0.05$m$^{-1}$. The trajectories are extended for both positive and negative $t.$}
\label{fig:q1vsq2}
\end{center}
\end{figure}

It is now reasonable to ask how the surface waves are affected by the vortex. They are governed by a KdV-type equation with variable coefficients and a source. This problem of course needs a proper numerical and analytical investigation, however our experience with similar equations (KdV with variable coefficients, reflecting variable bottom, \cite{CIT,CIMT,IMT}) indicates the following probable effects on an incoming KdV ``exact'' soliton. First, the soliton parameters such as amplitude and velocity could change, and second, a birth of a new soliton(s) may happen as well. These processes are usually accompanied by some energy loss in the form of radiation waves.

We mention also that the problem of surface-vortex interactions has been studied by different methods by Curtis and Kalisch in \cite{CKa}.

\section{Discussion}

We have explored the wave motion in a fluid domain with a free surface and a flat bottom of finite depth without any restrictions on the vorticity field.
As expected, the governing equations in this case can not be split into separate sets of equations for the free surface and the fluid volume.
Due to the nonconstant vorticity field there is always an interaction between the surface and body of the fluid and the equations are mathematically complicated.
Nevertheless, we have argued, that the equations could be in a relatively compact and manageable form by the means of the Dirichlet-Neumann operastors and Green's functions for the fluid domain. We have derived an expression for the Green functions which allows subsequent expansions for specific propagation regimes.

There are some challenging outstanding problems of purely mathematical nature, such as the symmetry of the Green function as well as the connection between the different representations of the Green function.

The complete Hamiltonian formulation of the problem is of fundamental importance for the water-wave theory and we intend to come back to this problem in
the near future.

We have illustrated the arising system of evolutionary equations by the example of a point vortex. Then this example has been further simplified for the Boussinesq and KdV regimes. The detailed further study of such systems will require numerical simulations and will be performed in follow-up publications.

\subsection*{Acknowledgments}

The authors would like to thank the Erwin Schr\"odinger International Institute for Mathematics and Physics
(ESI), Vienna (Austria) for the opportunity to participate in the Research in Teams project: Variational approaches to modelling geophysical waves and flows (July 1, 2022 — July 31, 2022) where a significant part of this work has been accomplished. RI acknowledges also partial funding from grant 21/FFP-A/9150 - Science Foundation Ireland.

%\subsection*{Data Availability Statement} Data sharing not applicable to this article as no datasets were generated or analysed during the current study.

\section*{Appendix: About the strength of the vortex}

Since the range of localization of the  surface variables $\mathbf{u}$, $\eta$ and the point  vortex is different, and a direct comparison would be misleading, let us now compare  the average values of the energy of both the surface wave and the surface motion due to the vortex over some interval, say the interval $[0,\lambda]$, where the wave amplitude, for example, has its maximum (since we, nevertheless consider a solitary wave).  The kinetic energy of the surface wave, $E_1,$ on average, is equal to its potential energy, $E_2$ (by absolute value)
\begin{align}\label{energies}
E_1:=\ \frac{1}{2\lambda}\int_0^\lambda \mathfrak{u}^2 dx, \quad E_2:=\  \frac{1}{2\lambda} \int_0^{\lambda}\frac{g}{h}\eta^2 dx,
\end{align}
respectively.
The energy of the motion, due to the vortex at the surface $z=0$, could be approximated by
\begin{align}\label{E1}
E_1\approx \frac{(\omega^*)^2}{2\lambda}\int_{-\infty}^{\infty} \big(\Gamma_x^2(\mathbf{x},\mathbf{q})+\Gamma_z^2(\mathbf{x},\mathbf{q})\big)\Big|_{z=0}dx.
\end{align}
For simplicity, we take
\begin{align}
\Gamma(\mathbf{x},\mathbf{z})=\frac{1}{4\pi}\log\big((x-q_1)^2+(z-q_2)^2\big),
\end{align}
and in view of the formula $\int_{-\infty}^{\infty}\frac{ dx}{x^2+a^2}=\frac{\pi}{|a|}$, we get for the right hand side in \eqref{E1} the expression
\begin{align}
 \frac{(\omega^*)^2}{8\lambda\pi^2}\int_{-\infty}^{\infty}
\frac{dx}{(x-q_1)^2+q_2^2}=  \frac{(\omega^*)^2}{8\lambda\pi| q_2|},\nonumber
\end{align}
hence,
\begin{align}\label{114}
E_1 \approx  \frac{(\omega^*)^2}{8\lambda\pi| q_2|}.
\end{align}
On the other hand, $E_2$ is proportional to $a^2$, indeed, if $\eta=a\cos(k(x-ct))$ for example, which is the expression in the linear case,
\begin{align}\label{115}
E_2\approx   \frac{1}{2\lambda} \int_0^{\lambda}\frac{g}{h} a^2\cos^2\left(\frac{2\pi (x-ct)}{\lambda}\right) dx= \frac{g a^2}{4 h},
\end{align}
where we used the formula $\int_0^{\lambda}\cos^2\left(\frac{2\pi (x-ct)}{\lambda}\right) dx=\frac{\lambda}{2}$. Comparing  \eqref{114} and  \eqref{115} yields
\begin{align}
(\omega^*)^2\sim \frac{2\pi \lambda g a^2 |q_2|}{ h}\implies \frac{(\omega^*)^2}{gh^3}\sim \frac{2\pi\lambda}{h}\left(\frac{a}{h}\right)^2\frac{|q_2|}{h}.
\end{align}
We introduce a parameter, related to the depth of the vortex, denoted by
\begin{align}
\mu:= \frac{|q_2|}{h}.
\end{align}
Thus, if $\mu\in \mathcal{O}(1)$, that is, the point vortex near to the bottom, then, for $\epsilon\sim \delta^2$,
\begin{align}\label{122}
\frac{\omega^*}{h\sqrt{gh}}\in \mathcal{O}(\delta^{\frac{3}{2}}).
\end{align}
If $\mu\in \mathcal{O}(\delta)$, that is, the point vortex is close to the free surface, then,
\begin{align}\label{123}
\frac{\omega^*}{h\sqrt{gh}}\in \mathcal{O}(\delta^2)=\mathcal{O}(\epsilon).
\end{align}

Therefore, the vortex strength which would make an impact on the equation, corresponding to the Boussinesq and KdV regimes, depends on the vortex depth and has approximate orders as in \eqref{122} and \eqref{123}.

\end{document}